\documentclass[11pt]{article}

\usepackage{amstext,amssymb,amsmath,amsbsy,bm}
\usepackage{hyperref}
\usepackage{amscd}
\usepackage{amsfonts}
\usepackage{indentfirst}
\usepackage{verbatim}
\usepackage{amsmath}
\usepackage{amsthm}
\usepackage{enumerate}
\usepackage{graphicx}
\usepackage{color}
\usepackage[OT1]{fontenc}
\usepackage[latin1]{inputenc}
\usepackage[english]{babel}
\usepackage{amssymb}
\usepackage{subfig}
\usepackage{algorithm}
\usepackage{algpseudocode}
\usepackage[shortlabels]{enumitem}

\newcommand{\R}{\mathbb{R}}
\newcommand{\ds}{\displaystyle}

\newcommand{\x}{{\bf x}}

\newcommand{\bu}{{\bf u}}
\newcommand{\bv}{{\bf v}}

\newcommand{\ii}{{\rm i}}
\newcommand{\ik}{\ii k}

\newcommand{\usc}{u_{\rm sc}}
\newcommand{\uinc}{u_{\rm inc}}

\newtheorem{Theorem}{Theorem}[section]
\newtheorem{Lemma}{Lemma}[section]

\newtheorem{Remark}{Remark}[section]

\newtheorem*{Assumption*}{Assumption}

\newtheorem{Problem}{Problem}[section]
\newtheorem*{Problem*}{Problem}
\numberwithin{equation}{section}

\usepackage{matlab-prettifier}
\usepackage{mathtools}

\begin{document}

\title{Inverse scattering without phase: Carleman convexification and phase retrieval via the Wentzel--Kramers--Brillouin approximation}

\author{
Thuy T. Le\thanks{%
Department of Mathematics, NC State University, Raleigh, NC 27695, USA, \texttt{tle9@ncsu.edu}.}
\and Phuong M. Nguyen\thanks{Department of Mathematics and Statistics, University of North Carolina at
Charlotte, Charlotte, NC, 28223, USA, \texttt{pnguye45@charlotte.edu}.
}
\and Loc H. Nguyen\thanks{Department of Mathematics and Statistics, University of North Carolina at
Charlotte, Charlotte, NC, 28223, USA, \texttt{loc.nguyen@charlotte.edu}.}
}

\date{}
\maketitle
\begin{abstract}
This paper addresses the challenging and interesting inverse problem of reconstructing the spatially varying dielectric constant of a medium from phaseless backscattering measurements generated by single-point illumination. The underlying mathematical model is governed by the three-dimensional Helmholtz equation, and the available data consist solely of the magnitude of the scattered wave field. To address the nonlinearity and servere ill-posedness of this phaseless inverse scattering problem, we introduce a robust, globally convergent numerical framework combining several key regularization strategies. Our method first employs a phase retrieval step based on the Wentzel--Kramers--Brillouin (WKB) ansatz, where the lost phase information is reconstructed by solving a nonlinear optimization problem. Subsequently, we implement a Fourier-based dimension reduction technique, transforming the original problem into a more stable system of elliptic equations with Cauchy boundary conditions. To solve this resulting system reliably, we apply the Carleman convexification approach, constructing a strictly convex weighted cost functional whose global minimizer provides an accurate approximation of the true solution. Numerical simulations using synthetic data with high noise levels demonstrate the effectiveness and robustness of the proposed method, confirming its capability to accurately recover both the geometric location and contrast of hidden scatterers.
\end{abstract}


\vspace{1em}

\section{Introduction}
\label{sec intr}

Inverse scattering problems, particularly those involving wave propagation through media with spatially varying dielectric properties, have broad applications in scientific and engineering disciplines. Consider a medium whose dielectric characteristics are represented by a smooth function \( c: \mathbb{R}^3 \to [1, \infty) \). Let \( \Omega = (-R, R)^3 \subset \mathbb{R}^3 \), with \( R > 0 \), denote our region of interest. Throughout this study, we assume that the dielectric constant satisfies
\begin{equation}
    c(\x) = \left\{
        \begin{array}{ll}
            \geq 1, & \x \in \Omega,\\[6pt]
            = 1, & \x \in \mathbb{R}^3 \setminus \Omega,
        \end{array}
    \right.
    \label{1.1}
\end{equation}
which implies a vacuum-like medium surrounding \(\Omega\), and regions of higher dielectric values within \(\Omega\) correspond to embedded scatterers. Practical examples of such scatterers include buried explosive devices, underground mineral deposits, tumors in biological tissues, defects within structural materials, and microscopic or nanoscopic inhomogeneities. These inverse scattering scenarios are fundamental to many areas, such as remote sensing, non-destructive evaluation, medical imaging, and materials science.

\begin{figure}[h!]
	\centering
	\subfloat[Schematic illustration of inverse scattering]{
		\includegraphics[width=.35\textwidth]{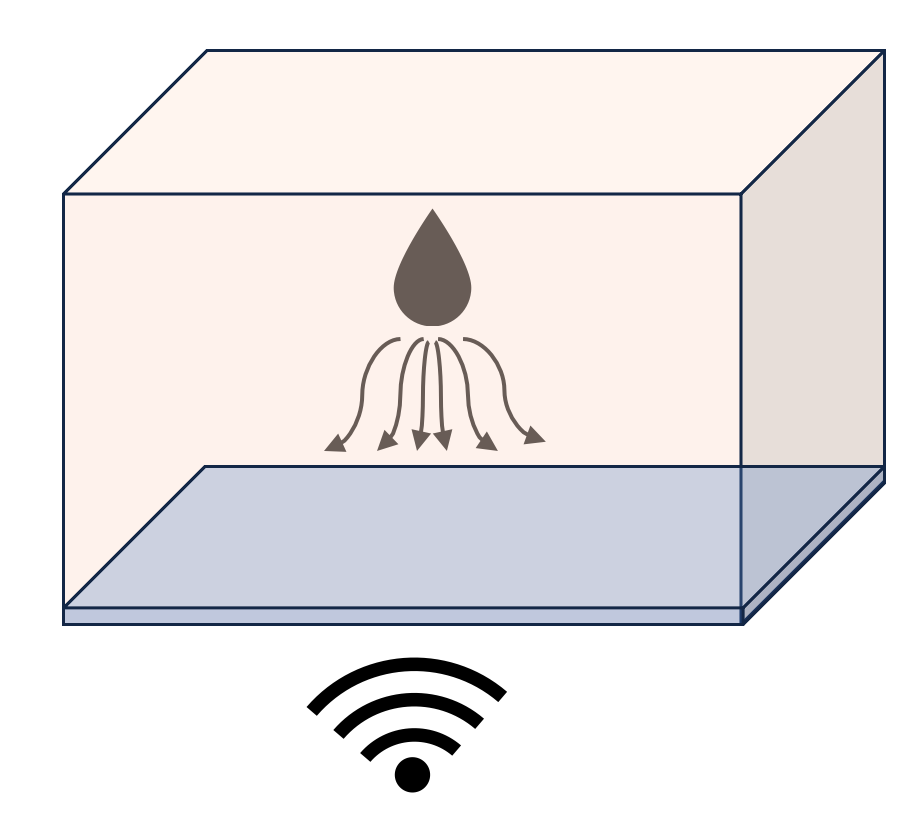}
		\put(-90,5){\tiny Incident wave}
		\put(-90,80){\tiny Scattered wave}
		\put(-95,105){\tiny Unknown object}
		\put(-140,135){\tiny Explored region $\Omega$}
		\put(-115,50){\tiny Measurement layer $\Gamma_L$}
	} \quad
	\subfloat[\label{fig1b}Remote sensing example]{
		\includegraphics[height=.35\textwidth,width=.35\textwidth]{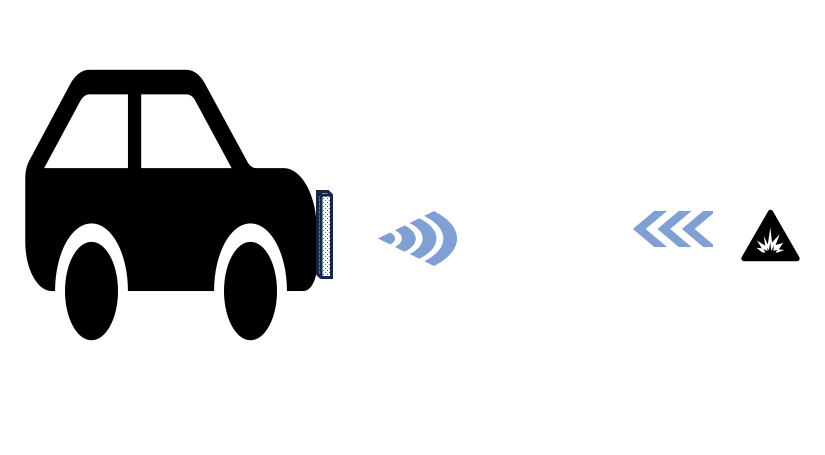}
		\put(-17,68){\tiny Unknown object}
		\put(-50,95){\tiny Backscattered wave}
		\put(-100,70){\tiny Incident wave}
		\put(-108,103){\tiny Detector array}
	}
	\caption{\label{diagram}
	(a) Illustration of the inverse scattering scenario, showing incident waves interacting with an unknown object in an inaccessible region \(\Omega\), producing scattered waves measured on the observation layer \(\Gamma_L\). (b) Real-world example demonstrating remote detection of hidden objects using reflected wave data collected by a detector array.}
\end{figure}
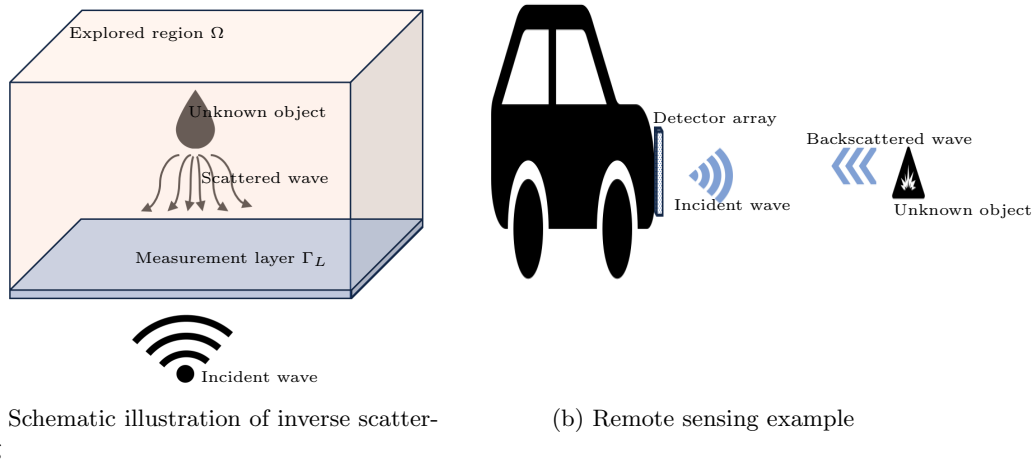

To remotely detect and characterize scatterers, one illuminates the domain \(\Omega\) with an incident wave originating from a single point source located at \(\x_0 = (0,0,-d) \notin \Omega\), where \( d > R \). This wave interacts with scatterers within \(\Omega\), creating scattered waves that radiate outward. We measure only the magnitude (intensity) of the resulting scattered wave field on a three-dimensional measurement layer
\[
	\Gamma_{L} = (-R, R)^2 \times (-R, -R + L),
\]
located immediately at the bottom of the domain \(\Omega\), with some thickness \( L > 0 \). Figure~\ref{diagram} provides both schematic and practical visualizations of this measurement scenario. The central goal of our phaseless inverse scattering problem is to reconstruct the dielectric function \( c(\x) \) using these intensity-only measurements, thereby revealing the position, shape, and dielectric contrasts of the hidden scatterers. 

Phaseless inverse scattering problems naturally arise in applications where the phase measurement is difficult or infeasible, such as in high-frequency regimes with rapid oscillations that preclude accurate phase detection. In contrast to conventional inverse scattering problems, where both amplitude and phase are typically available, phaseless measurements lead to greater mathematical complexity due to increased nonlinearity and limited data. As a consequence, it is generally necessary to collect intensity data over a three-dimensional measurement volume rather than on a simpler two-dimensional boundary surface.

The theoretical and practical significance of phaseless inverse scattering has long been recognized. A classical question posed by Chadan and Sabatier~\cite[Chapter 10]{Chadan1989} ``\emph{How can inverse scattering problems be solved using only magnitude data?}" has driven extensive research. Early theoretical investigations, such as uniqueness results for the one-dimensional Schr\"odinger equation~\cite{Klibanov:jmp1992,AktosunSacks:ip1998}, laid foundational understanding. Subsequently, these uniqueness results have been extended to multi-dimensional Schr\"odinger and Helmholtz equations under various assumptions~\cite{Klibanov:AppliedMathLetters2014,Klibanov:SiamJMath2014,KlibanovRomanov:ip2017,Klibanov:ipi2017,XuZhangZhang2018}. 

While analytic reconstruction methods have been proposed~\cite{KR:JIpp2015,KlibanovRomanov:ejmca2015,KlibanovRomanov:SIAMam2016,KlibanovRomanov:ip2016}, these methods often require measurements across a wide frequency spectrum, limiting their practicality. More computationally feasible numerical approaches for realistic frequency ranges have thus been developed, demonstrating successful reconstructions from simulated and experimental data~\cite{KlibanovKoshevNguyens:SJIS2018,KlibanovLiemLoc:SIAMjap2019, KlibanovLocKejia:apnum2016}. Further numerical methods utilizing Kirchhoff migration and Born approximations have been studied in~\cite{BardsleyVasquez:sjis2016,BardsleyVasquez:ip2016}, and shape reconstruction from phaseless data has also received considerable attention~\cite{AmmariChowZou:sjap2016,BaoZhang:ip2016,IvanyshynKress:jcp2011,IvanyshynKressSerranho:acm2010,LiLiuWang:jcp2014}.

In this paper, we investigate the phaseless inverse scattering problem under single-source illumination using multi-frequency data. While linearization-based methods can effectively localize scatterers, they often fail to accurately reconstruct high-contrast variations, see  e.g. \cite{KlibanovLocKejia:apnum2016}. To address this limitation, we propose a novel, globally convergent reconstruction method that combines the Wentzel--Kramers--Brillouin (WKB) ansatz, a Fourier filtering technique, and the Carleman convexification approach. By entirely avoiding linearization, our method significantly improves reconstruction accuracy and robustness.
Specifically, our proposed numerical framework consists of three interconnected components. First, we develop a robust phase retrieval procedure that leverages the WKB ansatz and nonlinear optimization to accurately recover complex wave fields from phaseless measurements. Next, we apply a Fourier-based frequency reduction technique, transforming the original inverse problem into a stable elliptic system with Cauchy boundary conditions. Finally, the Carleman convexification method is employed to construct a weighted cost functional, ensuring strict convexity and enabling stable global minimization. Numerical experiments clearly illustrate that our approach effectively recovers both the location and high contrast of scatterers, even in the presence of noise.

The remainder of the paper is organized as follows. In Section~\ref{sec_phase}, we outline the problem statement and our WKB-based phase retrieval procedure. Section~\ref{sec_reduceEqn} describes the frequency dimension reduction approach. Section~\ref{sec_convex} introduces the Carleman convexification method, including rigorous convexity analysis. Numerical simulations validating the proposed method are provided in Section~\ref{sec_num}, and concluding remarks are offered in Section~\ref{sec_concl}.

\section{Problem statement and phase reconstruction via the Wentzel--Kramers--Brillouin ansatz} \label{sec_phase}

Consider an inverse scattering scenario governed by the Helmholtz equation in three-dimensional space. Let $[\underline{k}, \overline{k}]$ denote an interval of wave numbers and $\mathbf{x}_0 = (0, 0, -d)$, with $d > R$, represent a fixed point source outside the domain $\Omega$. The incident wave generated by this source at wave number $k$ is given explicitly by
\begin{equation}
	u_{\mathrm{inc}}(\mathbf{x}, k) = \frac{e^{\mathrm{i} k |\mathbf{x} - \mathbf{x}_0|}}{4\pi |\mathbf{x} - \mathbf{x}_0|}, \quad (\mathbf{x}, k) \in \mathbb{R}^3 \times [\underline{k}, \overline{k}].
	\label{inc}
\end{equation}
When this incident wave interacts with the scatterers embedded in the domain $\Omega$, characterized by the spatially varying dielectric constant $c(\mathbf{x})$, it generates the total wave field $u(\mathbf{x}, k)$, which satisfies the following problem:
\begin{equation}
	\left\{
	\begin{array}{rcll}
		\Delta u + k^2 c(\mathbf{x}) u  &=& 0,  &\mathbf{x} \in \mathbb{R}^3,\\[6pt]
		\partial_{|\x|} u_{\mathrm{sc}} - \mathrm{i} k u_{\mathrm{sc}} &=& o(|\mathbf{x}|^{-1}), & |\mathbf{x}| \to \infty,
	\end{array}
	\right.	
	\label{Hel}
\end{equation}
where the scattered wave $u_{\mathrm{sc}}$ is defined as the difference between the total and incident waves, i.e.,
\begin{equation}
	u_{\mathrm{sc}}(\mathbf{x}, k) = u(\mathbf{x}, k) - u_{\mathrm{inc}}(\mathbf{x}, k), \quad (\mathbf{x}, k) \in \mathbb{R}^3 \times [\underline{k}, \overline{k}].
\end{equation}

In many practical situations, only the intensity (modulus) of the scattered wave can be measured. Thus, we consider the following inverse problem:

\begin{Problem}[Inverse scattering without phase information]	
	Given the phaseless measurements
	\begin{equation}
		f(\mathbf{x}, k) = |u(\mathbf{x}, k)|, \quad (\mathbf{x}, k) \in \Gamma_L \times [\underline{k}, \overline{k}],
		\label{fmea}
	\end{equation}
	reconstruct the dielectric constant $c(\mathbf{x})$ within the domain $\Omega$.
	\label{cip}
\end{Problem}

Addressing this phaseless inverse scattering problem presents significant challenges due to its inherent nonlinearity and ill-posedness. To reconstruct the lost phase information, we utilize the classical Wentzel--Kramers--Brillouin (WKB) ansatz, which approximates the wave field:
\begin{equation}
    u(\mathbf{x}, k) = A(\mathbf{x}) e^{i k \tau(\mathbf{x})} + \mathcal{O}(1/k) 
    \quad \text{as } k \to \infty.
    \label{2.1}
\end{equation}

This asymptotic form is well-established in the literature. For the stationary Schr\"odinger equation, its validity under the Born approximation has been demonstrated in foundational studies such as \cite{Critchfield1964, Hecht1957, Miller1953}, and analogous reasoning extends naturally to the Helmholtz equation. Beyond the Born approximation, the first rigorous justification was presented in \cite[Theorem 17]{Vainberg:rms1975}, with a more refined analysis appearing much later in \cite{KlibanovRomanov:SIAMam2016}. For clarity, we summarize below the key conditions provided in \cite{KlibanovRomanov:SIAMam2016} under which the ansatz is valid:

\begin{enumerate}
    \item The coefficient \( c(\mathbf{x}) \) is assumed to belong to the class \( C^{15} \).
    \item For every \( \mathbf{x} \in \Omega \), there exists a unique geodesic line connecting  the source \( \mathbf{x}_0 \) and \( \mathbf{x} \) with respect to the Riemannian metric \( \sqrt{c(\mathbf{x})} |d\mathbf{x}| \), where \( |d\mathbf{x}| = \sqrt{dx^2 + dy^2 + dz^2} \).
\end{enumerate}

Although these conditions rigorously justify the ansatz \eqref{2.1}, they may not be practical to verify for finite \( k \) or in real-world applications. Therefore, we adopt the ansatz heuristically as a physically informed starting point for phase reconstruction.
To estimate \( \tau \), we differentiate the ansatz, neglecting the \( \mathcal{O}(1/k) \) term:
\begin{align*}
	\nabla u(\x, k) &= \left(\nabla A(\x) + i k \nabla \tau(\x)\right) e^{i k \tau(\x)},
\\
	\Delta u(\x, k) &= \left(
		\Delta A(\x) + 2i k \nabla A(\x) \cdot \nabla \tau(\x) + i k A(\x)\Delta \tau(\x) - k^2 A(\x) |\nabla \tau(\x)|^2
	\right) e^{i k \tau(\x)}.
\end{align*}
Substituting into the Helmholtz equation, we obtain:
\begin{multline}
	\Delta u(\x, k) + k^2 c(\x) u(\x, k) = \Big(
		\Delta A(\x) + 2i k \nabla A(\x) \cdot \nabla \tau(\x) 
		\\
		+ i k A(\x)\Delta \tau(\x) - k^2 A(\x) |\nabla \tau|^2 + k^2 c(\x) A(\x)
	\Big) e^{i k \tau(\x)} = 0,
	\label{2.2}
\end{multline}
for all \( (\x, k) \in \Gamma_L \times [\underline k, \overline k] \). In the high-frequency limit, the leading-order term,
\[
    -k^2 A(\x) |\nabla \tau|^2 + k^2 c(\x) A(\x),
\]
must vanish, which yields the eikonal equation:
\begin{equation}
    |\nabla \tau|^2 = c(\x) = 1 \quad \text{for } \x \in \Gamma_L,
    \label{eikonal}
\end{equation}
using the assumption \( c = 1 \) in \( (\mathbb{R}^3 \setminus \Omega) \cup \Gamma_L \).

\begin{Remark}
	Alternatively, the eikonal equation can be heuristically derived by viewing \eqref{2.2} as a quadratic equation in \( k \). Enforcing that the equation holds for all \( k \) implies that the coefficient of the dominant term (and the other two lower order terms) must vanish, once again yielding \eqref{eikonal}.
\end{Remark}

Since \( \tau \) denotes the travel time from \( \x_0 \) to \( \x \), the natural choice of solution to \eqref{eikonal} is:
\begin{equation}
    \tau(\x) = |\x - \x_0|  
    \quad
    \text{for all } \x \in \Gamma_L.
    \label{2.2222}
\end{equation}
Using this phase, we define the initial guess:
\begin{equation}
    u^{(0)}(\x, k) = f(\x, k) e^{i k |\x - \x_0|}
    \quad \text{for all } \x \in \Gamma_L, \quad k \in [\underline k, \overline k],
    \label{2.5555}
\end{equation}
as an approximation to the true wave field \( u(\x, k) \). 

The ansatz \eqref{2.1} provides a mathematically grounded approximation of wave fields in the high-frequency regime. However, the constructed initial guess \( u^{(0)}(\x, k) \) may not exactly satisfy the Helmholtz equation, nor perfectly match the observed modulus data due to noise. To refine this estimate, we recover an improved approximation \( u_{\rm phase}(\x, k) \) by minimizing its deviation from both the Helmholtz operator and the intensity constraint. Specifically, for each \( k \in [\underline k, \overline k] \), we solve the following variational problem:

\begin{equation}
    J_k(v) 
    = \big\|\Delta v + k^2 v\big\|_{L^2(\Gamma_L)}^2 + \big\||v|^2 - f^2(\x, k)\big\|_{L^2(\Gamma_L)}^2,
    \label{2.5}
\end{equation}
where the minimization is over \( v \in H^2(\Gamma_L) \). Among the local minimizers, we select the one closest to \( u^{(0)}(\x, k) \), and denote it by \( u_{\rm phase}(\x, k) \).

The functional \( J_k(v) \) in \eqref{2.5} captures two complementary objectives: adherence to the PDE model in \( \Gamma_L \), and conformity with the observed intensity. The first term, \( \|\Delta v + k^2 v\|_{L^2(\Gamma_L)}^2 \), penalizes deviation from the Helmholtz equation, while the second term, \( \||v|^2 - f^2(\x, k)\|_{L^2(\Gamma_L)}^2 \), ensures alignment with the modulus data. This dual-objective formulation is standard and widely used in scientific and engineering communities. 
Restricting to the Sobolev space \( H^2(\Gamma_L) \) ensures well-posedness and the necessary regularity to evaluate Laplacians and traces.
This framework is particularly advantageous for mitigating the ill-posedness introduced by noise in the measured data.
 Initializing the optimization at \( u^{(0)} \), which incorporates travel-time-informed phase information, increases the likelihood of convergence to a physically relevant solution in this non-convex landscape. 
 
Observe that \( \Gamma_L \subset \Omega \), and its bottom boundary coincides with:
\begin{equation}
    \Gamma = (-R, R)^2 \times \{ z = -R \},
\end{equation}
which corresponds to the lower face of \( \partial \Omega \). Thus, both \( u_{\rm phase} \) and its normal derivative \( \partial_z u_{\rm phase} \) can be extracted on \( \Gamma \), enabling the reformulation of the original phaseless inverse problem as follows:

\begin{Problem}[The inverse scattering problem with phase information]
Given the functions
\begin{equation}
	g(\x, k) = u(\x, k) \quad 
	\text{and} \quad
	h(\x, k) = \partial_z u(\x, k)
	\label{2.6}
\end{equation}
for all \( \x \in \Gamma \times [\underline{k}, \overline{k}] \),
reconstruct \( c(\x) \) for all \( \x \in \Omega \).
\label{p1_tran}
\end{Problem}

This formulation pertains to the phased inverse scattering problem, where both the modulus and phase of the backscattered wave are known on the measurement surface. Our research group has developed two distinct approaches, both leveraging Carleman estimates, to address this problem:  
\begin{enumerate}
	\item The Carleman convexification method \cite{VoKlibanovNguyen:IP2020, Khoaelal:IPSE2021, KhoaKlibanovLoc:SIAMImaging2020, LeNguyen:JSC2022};
	\item The Carleman contraction mapping method \cite{NguyenNguyenVu}.
\end{enumerate}
In this work, we adopt the first approach to solve the problem numerically. The primary motivation for this choice is that the Carleman convexification method has been validated with both simulated and experimental datasets. The contraction mapping method will be applied in future research.

\section{The frequency dimension reduction model} \label{sec_reduceEqn}

Following the phase retrieval procedure described in Section~\ref{sec_phase}, Problem~\ref{cip} is reduced to reconstructing the coefficient \( c(\x) \) from boundary measurements of \( u(\x, k) \) and \( \partial_z u(\x, k) \) on \( \Gamma \times [\underline{k}, \overline{k}] \), as defined in Problem~\ref{p1_tran}. This inverse scattering problem is highly ill-posed, meaning that even minor perturbations in the input data, such as measurement noise, can result in substantial reconstruction errors. To mitigate this instability, we employ a frequency filtering technique using Fourier truncation, which suppresses high-frequency oscillations in the data. This process yields a system of elliptic equations with Cauchy boundary conditions, which is notably more stable and amenable to numerical computation. The formulation of this system also naturally aligns Problem~\ref{p1_tran} with the framework of the Carleman convexification method.

\subsection{The logarithmic transformation}

For the reader's convenience, we provide an overview of the Fourier truncation procedure in this section. The process begins with an algebraic transformation, following the algorithmic approach developed in \cite{Khoaelal:IPSE2021, KhoaKlibanovLoc:SIAMImaging2020, KlibanovLiemLocHui:IPI2018, LeNguyen:JSC2022}. Specifically, for each \( (\x, k) \in \Omega \times [\underline k, \overline k] \), we define the logarithmic transformation
\begin{equation}
    v(\x, k) = \frac{1}{k^2} \log \frac{u(\x, k)}{\uinc(\x, k)}.
    \label{change_variable}
\end{equation}
\begin{Remark}
Although taking the logarithm of the complex-valued function \( \frac{u(\x, k)}{\uinc(\x, k)} \) in \eqref{change_variable} may initially appear problematic, it is well-defined based on the definition of the complex logarithm presented in \cite[Section 4.2]{KlibanovLiemLocHui:IPI2018}. The WKB ansatz \eqref{2.1} further supports this definition by ensuring that \( \frac{u(\x, k)}{\uinc(\x, k)} \) remains nonzero for all \( (\x, k) \in \Omega \times [\underline{k}, \overline{k}] \). In numerical implementations, evaluating this logarithm is straightforward and poses no practical difficulties.
\end{Remark}

By standard rules in differentiation, we have for all $(\x, k) \in \Omega \times [\underline k, \overline k]$
\begin{equation*}
	\nabla v(\x, k) = \frac{1}{k^2} \Big[
		\frac{\nabla u(\x, k)}{u(\x, k)} - \frac{\nabla \uinc(\x, k)}{\uinc(\x, k)}
	\Big],
\end{equation*}
and
\begin{align}
	\Delta v(\x, k) &= \frac{1}{k^2} \Big[
		\frac{\Delta u(\x, k)}{u(\x, k)} - \Big(\frac{\nabla u(\x, k)}{u(\x, k)}\Big)^2 -  \frac{\Delta \uinc(\x, k)}{\uinc(\x, k)} + \Big(\frac{\nabla \uinc(\x, k)}{\uinc(\x, k)}\Big)^2
	\Big] \notag
	\\
	&= \frac{1}{k^2} \Big[
		k^2(1 - c(\x)) - \Big(
			\frac{\nabla u(\x, k)}{u(\x, k)} - \frac{\nabla \uinc(\x, k)}{\uinc(\x, k)}
		\Big) \cdot\Big(
			\frac{\nabla u(\x, k)}{u(\x, k)} + \frac{\nabla \uinc(\x, k)}{\uinc(\x, k)} 
		\Big)\Big] \notag
		\\
	&= 1 - c(\x) -  \nabla v(\x, k) \cdot \Big(
			k^2\nabla v(\x, k) + \frac{2\nabla \uinc(\x, k)}{\uinc(\x, k)}
		\Big).	 \label{3.2}
\end{align}
By a direct algebra, using the explicit formula of $\uinc$ in \eqref{inc} gives 
\begin{equation}
	\frac{\nabla \uinc(\x, k)}{\uinc(\x, k)}= \Big(\ik - \frac{1}{|\x - \x_0|}\Big) \frac{\x - \x_0}{|\x - \x_0|},
	\quad (\x, k) \in \Omega \times [\underline k, \overline k].
	\label{3.3}
\end{equation}
Combining \eqref{3.2} and \eqref{3.3} gives 
\begin{equation}
	\Delta v(\x, k) + k^2 [\nabla v(\x, k)]^2 + 2 \Big(\ik - \frac{1}{|\x - \x_0|}\Big) \nabla v(\x, k) \cdot \frac{\x - \x_0}{|\x - \x_0|}  = 1 - c(\x)
	\label{3.4}
\end{equation}
for all $(\x, k) \in \Omega \times [\underline k, \overline k]$.
To eliminate the unknown $c$, we differentiate \eqref{3.4} with respect to $k$ to obtain
\begin{multline}
	\Delta \partial_k v(\x, k) + 2k [\nabla v(\x, k)]^2 + 2k^2 \nabla v(\x, k) \cdot \nabla \partial_k v(\x, k) 
	\\
	+ 2 \Big(\ik - \frac{1}{|\x - \x_0|}\Big)  \nabla \partial_k v(\x, k) \cdot \frac{\x - \x_0}{|\x - \x_0|}
	+ 2\ii\nabla v(\x, k) \cdot \frac{\x - \x_0}{|\x - \x_0|} = 0
	\label{3.5}
\end{multline}
for $\x \in \Omega,$ $k \in [\underline k, \overline k]$.

Solving \eqref{3.5} is challenging because it does not take the form of a standard partial differential equation, and a theoretical framework for it has not yet been developed. As a result, we approach this problem using a frequency dimension reduction technique.

\subsection{The frequency dimension reduction using Fourier expansion}

We next apply a ``Fourier filter" to eliminate the high-frequency oscillatory components of the function \( v \). This process involves truncating the Fourier expansion of \( v \) using the polynomial-exponential basis \( \{\Psi_n\}_{n \geq 1} \) of \( L^2(\underline{k}, \overline{k}) \). This basis was originally constructed in \cite{Klibanov:jiip2017} via the Gram-Schmidt orthonormalization of the complete system \( \{\phi_n(k) = k^{n-1} e^k\}_{n \geq 1} \) in \( [\underline k, \overline k]\). A higher-dimensional extension of this basis was later developed in \cite{NguyenLeNguyenKlibanov:2023}. The rationale behind this choice will be further discussed in Remark~\ref{rem 21}.

For $\x \in \Omega$, the Fourier expansion of $v$ is approximated as
\begin{equation}
	v(\x, k) = \sum_{n = 1}^\infty v_n(\x) \Psi_n(k) \approx \sum_{n = 1}^N v_n(\x) \Psi_n(k)
	\label{3.6}
\end{equation} 
where the cutoff number $N$ will be chosen later by a trial-error procedure, and the Fourier coefficient $v_n$ is given by
\begin{equation}
	v_n(\x) = \int_{\underline k}^{\overline k} v(\x, k) \Psi_n(k) dk.
	\label{3.6666}
\end{equation}

\begin{Remark}
The truncation in \eqref{3.6} acts as a filtering step, effectively regularizing any noise in the measured data by completely removing the high-frequency oscillatory components of $v$. Additionally, it significantly reduces computational costs by eliminating the frequency dimension. Specifically, instead of computing the function $v: \Omega \times [\underline{k}, \overline{k}] \to \mathbb{C}$, which involves $3+1$ dimensions, we compute a finite number $N$ of Fourier coefficients $v_{n}: \Omega \to \mathbb{C}$, reducing the problem to only three dimensions.
\label{rem31}
\end{Remark}

Plugging the approximation \eqref{3.6} into \eqref{3.5} gives
\begin{multline}
	 \sum_{n = 1}^N \Delta v_n(\x) \Psi_n'(k) 
	 + 2k \Big[\sum_{n = 1}^N \nabla v_n(\x) \Psi_n(k)\Big]^2 
	 + 2k^2 \sum_{n = 1}^N \nabla v_n(\x) \Psi_n(k) \cdot \sum_{n = 1}^N \nabla v_n(\x) \Psi_n'(k)
	\\
	+ 2 \Big(\ik - \frac{1}{|\x - \x_0|}\Big)  \sum_{n = 1}^N \nabla v_n(\x) \Psi_n'(k) \cdot \frac{\x - \x_0}{|\x - \x_0|}
	+ 2\ii\sum_{n = 1}^N \nabla v_n(\x) \Psi_n(k) \cdot \frac{\x - \x_0}{|\x - \x_0|} = 0
	\label{3.8}
\end{multline}
for all $(\x, k) \in \Omega \times [\underline k, \overline k].$
For each $m \in \{1, 2, \dots, N\}$, multiplying $\Psi_m(k)$ to both sides of \eqref{3.8} and integrating the resulting equation, we obtain
\begin{equation}
	 \sum_{n = 1}^N s_{mn} \Delta v_n(\x) 
	 +  \sum_{ n = 1}^N\sum_{ l = 1}^N a_{mnl}\nabla v_n(\x) \cdot \nabla v_l(\x) 	 	
	+   \sum_{n = 1}^N  {\bf b}_{mn} \cdot \nabla v_n(\x) 	 = 0
	\label{3.9}
\end{equation} for all $\x \in \Omega,$
where
\begin{align}
	s_{mn} &= \int_{\underline k}^{\overline k} \Psi_n'(k)  \Psi_m(k) dk, \label{smn}\\
	a_{mnl} &= \int_{\underline k}^{\overline k} 2k\Psi_n(k) \Psi_m(k)\big(\Psi_l(k) + k\Psi_l'(k) \big) dk,\\
	{\bf b}_{mn}(\x) &=  \frac{\x - \x_0}{|\x - \x_0|} \int_{\underline k}^{\overline k} 2 \Big[\Big(\ik - \frac{1}{|\x - \x_0|}\Big) \Psi_n'(k) + \ii \Psi_n(k)\Big]\Psi_m(k) dk \notag\\
	&=\frac{2\ii(\x - \x_0)}{|\x - \x_0|} \int_{\underline k}^{\overline k} k \Psi_n'(k) \Psi_m(k) dk
	- \frac{2(\x - \x_0) s_{mn}}{|\x - \x_0|^2} + \frac{2\ii(\x - \x_0)}{|\x - \x_0|}\delta(m, n).
\end{align}
Here,  $\delta(m, n)$ is the  Kronecker delta
\[
\delta(m, n) =
\begin{cases} 
1 & \text{if } m = n, \\
0 & \text{if } m \neq n.
\end{cases}
\]

\subsection{Data complementation}

The values of \( \bv(\x) = \left[\begin{array}{ccc} v_1(\x), \dots, v_N(\x) \end{array}\right]^{\rm T} \) for \( \x \in \Gamma \) can be explicitly computed as follows. From equations \eqref{2.6}, \eqref{change_variable}, and \eqref{3.6666}, we derive  
\begin{equation}
    g_m(\x) := v_m(\x) = \int_{\underline{k}}^{\overline{k}} \frac{\Psi_m(k)}{k^2} \log \frac{g(\x, k)}{\uinc(\x, k)} \, dk,
    \quad m \in \{1, 2, \dots, N\},
    \label{3.10}
\end{equation}  
and  
\begin{equation}
    h_m(\x) := \partial_z v_m(\x) = \int_{\underline{k}}^{\overline{k}} \frac{\Psi_m(k)}{k^2} \left[
        \frac{h(\x)}{g(\x)} - \left( \ik - \frac{1}{|\x - \x_0|} \right) \frac{z + d}{|\x - \x_0|}
    \right] dk,
    \label{3.1010}
\end{equation}
for all \( \x = (x, y, z) \in \Gamma \) and \( m \in \{1, 2, \dots, N\} \), where \( \x_0 = (0, 0, -d) \). While these expressions allow for the direct computation of \( \bv \) and its normal derivative on \( \Gamma \), they are not sufficient to determine \( \bv \) throughout the entire domain \( \Omega \) because
 solving second-order equations such as \eqref{3.9} requires knowledge of the solution on a surface enclosing \(\Omega\). However, the data in \eqref{3.10} and \eqref{3.1010} is well-defined only on \(\Gamma\), which is located at the bottom \(\Omega\). This indicates that the information of \(\{v_m\}_{m = 1}^N\) on \(\partial \Omega \setminus \Gamma\) is crucial but missing.  
To achieve a stable computation of \(\{v_m\}_{m = 1}^N\), it is necessary to compensate for this missing data. Since the scattering wave weakens when the receivers are far from the source, we can approximate the scattered wave \(\usc(\x, k)\) as zero on \(\partial \Omega \setminus \Gamma\). Consequently, \(v(\x, k)\) vanishes on \(\partial \Omega \setminus \Gamma \times [\underline{k}, \overline{k}]\). As a result, for \(m = 1, 2, \dots, N\), we obtain  
\begin{equation}
	v_m(\x) = 0 \quad \mbox{for all } \x \in \partial \Omega \setminus \Gamma.
	\label{data_supplement}
\end{equation}

\begin{Remark}
	The data supplementation described above is not rigorous but rather an approximation based on the fact that the scattering wave on \(\Gamma\) is stronger than that on \(\partial \Omega \setminus \Gamma\). This observation holds because the source location \(\x_0 = (0, 0, -d)\) is closer to \(\Gamma\) compared to \(\partial \Omega \setminus \Gamma\). This supplementation strategy has been successfully applied in our previous works on inverse scattering problems, such as \cite{VoKlibanovNguyen:IP2020, KhoaKlibanovLoc:SIAMImaging2020,  LeNguyen:JSC2022}. Therefore, we continue to employ this approach.
\end{Remark}

Combining \eqref{3.9}, \eqref{3.10}, and \eqref{data_supplement} gives
\begin{equation}
	\left\{
		\begin{array}{ll}
			\ds\sum_{n = 1}^N s_{mn} \Delta v_n(\x) 
	 +  \sum_{ n = 1}^N\sum_{ l = 1}^N a_{mnl}\nabla v_n(\x) \cdot \nabla v_l(\x) 	 	
	+   \sum_{n = 1}^N  {\bf b}_{mn} \cdot \nabla v_n(\x) 	 = 0 &\x \in \Omega,
	\\
	v_m(\x) = g_m(\x) &\x \in \Gamma,	
	\\
		\partial_z v_m(\x) = h_m(\x) &\x \in \Gamma,	
	\\
	v_m(\x) = 0 &\x \in \partial \Omega \setminus \Gamma.
		\end{array}
	\right.
	\label{3.11}
\end{equation}


\begin{Remark}
The system composed of all \( N \) equations in \eqref{3.11}, for \( m \in \{1, \dots, N\} \), provides an approximate model for addressing the inverse scattering problem. This approximation arises from the truncation of high-frequency components in \eqref{3.6} and the data complementation procedure described in \eqref{data_supplement}. Although these steps introduce inaccuracy, we view them as a deliberate and necessary regularization strategy. This trade-off substantially mitigates the inherent ill-posedness of the inverse problem, enabling a more stable and feasible numerical solution.
\end{Remark}

\section{The Carleman convexification method}\label{sec_convex}

To solve \eqref{3.11}, we employed the Carleman convexification method originally developed in \cite{KlibanovIoussoupova:SMA1995} and further advanced in subsequent works \cite{VoKlibanovNguyen:IP2020, Khoaelal:IPSE2021, KhoaKlibanovLoc:SIAMImaging2020, LeNguyen:JSC2022}.

Recall the matrix \( S = (s_{mn})_{m,n = 1}^N \) as the \( N \times N \) matrix whose \( (m,n) \)-entry is given by \eqref{smn}. From \cite{Klibanov:jiip2017}, it is known that \( S \) is invertible. We denote its inverse by \( S^{-1} = (\widetilde s_{mn})_{m, n = 1}^N \).  
Define 
\[
\widetilde{\bf b}_{mn} = \sum_{i = 1}^N\widetilde s_{mi} {\bf b}_{in}, \quad
\widetilde a_{mnl} = \sum_{i = 1}^N\widetilde s_{mi} a_{inl}, 
\quad \text{for } m, n, l \in \{1, \dots, N\}.
\]
Using \eqref{3.9}, we obtain the equation  
\begin{equation}
    \Delta  v_m(\x) 
    +  \sum_{n = 1}^N \nabla  v_n(\x)  \cdot \widetilde{\bf b}_{mn}		 
    + \sum_{n = 1}^N \sum_{l = 1}^N \widetilde  a_{mnl}\nabla v_n(\x)  \cdot  \nabla v_l(\x) 
    = 0,
    \label{4.1}
\end{equation}
for \( \x \in \Omega \).  
Equation \eqref{4.1}, combined with \eqref{3.10}, forms an elliptic system:  
\begin{equation}
    \left\{
        \begin{array}{ll}
            \ds\Delta  v_m(\x) 
            +  \sum_{n = 1}^N \nabla  v_n(\x)  \cdot \widetilde{\bf b}_{mn}		 
            + \sum_{n = 1}^N \sum_{l = 1}^N \widetilde  a_{mnl}\nabla v_n(\x)  \cdot  \nabla v_l(\x) 
            = 0, & \x \in \Omega,\\
            v_m(\x) = g_m(\x), & \x \in \Gamma,\\
            \partial_z v_m(\x) = h_m(\x), & \x \in \Gamma,\\
            v_m(\x) = 0 &\x \in \partial \Omega \setminus \Gamma.
        \end{array}
    \right.
    \label{4.2}
\end{equation}

\begin{Remark}
The system defined in \eqref{4.2} for \( m \in \{1, 2, \dots, N\} \) plays a key role to place Problem \ref{p1_tran} into the framework of the convexification method, and the choice of the basis \( \{\Psi_n\}_{n \geq 1} \) is crucial to its formulation. For a given \( N \in \mathbb{N} \), recall the matrix \( S \in \mathbb{R}^{N \times N} \) whose \( (m,n) \)-entry is given by \eqref{smn}. It was shown in \cite{Klibanov:jiip2017} that:
\begin{enumerate}
    \item The matrix \( S \) is invertible;
    \item For all \( n \geq 1 \), the derivative \( \Psi_n' \) is not identically zero on the interval \( [\underline k, \overline k] \).
\end{enumerate}

The first property is essential, as the invertibility of \( S \) is necessary to define the coefficients in \eqref{4.2}; without knowledge of \( S^{-1} \), the formulation would be incomplete. The second property is equally significant. If \( \Psi_n' = 0 \) for some \( n \), then important information about \( \Delta v_n \) in the first term of \eqref{3.11} would be lost. For example, consider replacing the exponential-polynomial basis with more commonly used alternatives such as the Legendre polynomials or the trigonometric basis \( \{\phi_n\}_{n \geq 1} \). In these cases, \( \phi_1 \) is a constant function, and thus \( \phi_1' \equiv 0 \). As a result, the principal term \( \Delta v_1 \) would be absent from \eqref{3.11}, leading to potentially large errors in computing \( v_1 \), and consequently affecting the accuracy of the reconstructed function \( v(\x, k) \) via \eqref{3.6}. This issue is particularly critical because the first term in the truncated series often contributes significantly to the overall reconstruction.
\label{rem 21}
\end{Remark}

Let \( s \geq 4\) be an integer. Note that $s$ is set to be larger than or equal to $4$ to ensure that \( H^s(\Omega) \) is continuously embedded into \( C^2(\overline \Omega) \). For theoretical purposes, we assume that the target function \( c \) is sufficiently smooth so that the true solution to the forward problem \( u \) satisfies  
\[
\|u(\cdot, k)\|_{H^s(\Omega)} < \infty, \quad \text{uniformly for } k \in [\underline k, \overline k].
\]
Additionally, we assume that \( |u(\x, k)| \) is uniformly bounded below by a positive constant \( u_0 \). Due to the change of variables in \eqref{change_variable}, the function \( v(\x, k) \) belongs to \( H^s(\Omega) \), with its \( H^s \) norm uniformly bounded for all \( k \in [\underline k, \overline k] \).  
Thus, there exists a constant \( M \), depending on the upper bound of \( \|u(\cdot, k)\|_{H^s(\Omega)} \) for \( k \in [\underline k, \overline k] \) and the set \( \{\Psi_n\}_{n = 1}^N \), such that  
\[
\|{\bf v}^*\|_{H^s(\Omega)^N} < M,
\]
where \( {\bf v}^* \) is the true solution to \eqref{4.2}.  
We seek the ``best fit" solution to \eqref{4.2} in the set of admissible solutions
\begin{equation}
	B(0, M) = \big\{
			\bm{\varphi}
			 \in H^s(\Omega)^N: \|\bm{\varphi}\|_{H^s(\Omega)^N} < M
	\big\},
\end{equation}
which is the ball centered at the origin with radius $M$.
For positive numbers $\lambda$ and $\epsilon$, 
define the Carleman weighted functional $J_{\lambda, \epsilon}: B(0, M) \to \R$
\begin{multline}
	J_{\lambda, \epsilon}(\bm{\varphi}) = \sum_{m = 1}^N \Big[\int_{\Omega} e^{2\lambda (z - r)^2}\Big|
		\Delta  \varphi_m(\x) 
            +  \sum_{n = 1}^N \nabla  \varphi_n(\x)  \cdot \widetilde{\bf b}_{mn}		 
            + \sum_{n = 1}^N \sum_{l = 1}^N \widetilde  a_{mnl}\nabla \varphi_n(\x)  \cdot  \nabla \varphi_l(\x)
	\Big|^2d\x
	\\
	+ \lambda^4 e^{2\lambda (R+r)^2}\int_{\Gamma} \big( |\varphi_m - g_m|^2 + |\partial_z \varphi_m - h_m|^2\big)d\sigma(\x)
		+ \lambda^4 \int_{\partial \Omega \setminus \Gamma} e^{2\lambda(z- r)^2} |\varphi_m|^2 d\sigma(\x)
		\\
	+ \epsilon \|\varphi_m\|_{H^s(\Omega)}^2\Big]
	\label{4.4}
\end{multline}
for all $\bm{\varphi} \in B(0, M)$,
where $r > R$ is a fixed number and $\lambda$ is a Carleman parameter. 
\begin{Remark}
 In the absence of the Carleman weight function \( e^{2\lambda (z - r)^2} \), the functional above reduces to the standard least-squares mismatch functional associated with \eqref{4.2}. However, in this unweighted form, this cost functional may possess multiple local minima, making its minimization particularly difficult. Conventional optimization techniques, such as the gradient descent method, are prone to becoming trapped in local minima that are far from the global solution. The core idea of the convexification method is to incorporate the Carleman weight function, which transforms the cost functional into a globally convex form, thereby facilitating the application of standard optimization methods and improving convergence to the true solution.
 \end{Remark}

The theoretical foundation of the convexification method is built upon a rigorous theorem that guarantees \( J_{\lambda, \epsilon} \) is strictly convex within the ball \( B(0, M) \), and that its global minimizer in this set provides a reliable approximation to the true solution \( \bv^* \) of \eqref{4.2}. This convexification theorem is fundamentally derived from the following Carleman estimate.

\begin{Lemma}[Carleman Estimate]
There exist constants $\lambda_0 = \lambda_0(\Omega, r) \geq 1$ and $C = C(\Omega, r) > 0$ such that 
\begin{multline} \int_{\Omega} e^{2\lambda (z - r)^2} |\Delta \varphi|^2,d\x \geq C \int_{\Omega} e^{2\lambda (z - r)^2}(\lambda^3|\varphi|^2 + \lambda|\nabla \varphi|^2)d\x
\\
-C\lambda^3\int_{\partial \Omega} e^{2\lambda (z - r)^2}  |\varphi|^2  d\sigma(\x)
-C\lambda \int_\Gamma e^{2\lambda(z - r)^2} |\nabla \varphi|^2d\sigma(\x).
\label{4.55}
\end{multline} 
for all functions $\varphi \in C^2(\Omega)$.
\label{Lemma_Car} 
\end{Lemma}
The Carleman estimate stated in Lemma~\ref{Lemma_Car} is structurally similar to the one presented in \cite[Theorem 4.1]{KlibanovLiZhang:SIAM2019}. The key distinction lies in the boundary conditions imposed on the function \( \varphi \). Specifically, \cite[Theorem 4.1]{KlibanovLiZhang:SIAM2019} assumes that \( \varphi|_{\partial \Omega} = 0 \) and \( \partial_z \varphi|_{\Gamma} = 0 \), whereas Lemma~\ref{Lemma_Car} omits the these conditions. To compensate for this relaxation, two negative terms are added to the right-hand side of the Carleman estimate in \eqref{4.55}.
Despite this difference, the proof technique used in \cite{KlibanovLiZhang:SIAM2019} remains applicable to Lemma~\ref{Lemma_Car} with only minor modifications. The main adjustment involves deferring the integration process. Instead of integrating mid-proof, as done in \cite[Theorem 4.1]{KlibanovLiZhang:SIAM2019}, we maintain the Carleman estimate locally at each point in \( \Omega \), deriving a pointwise estimate first and performing the integration only at the final stage. This postponed integration strategy aligns with the approach in \cite[Theorem 3.1 and Corollary 3.2]{LeLeNguyen:2024}. 
Given the similarity in methodology and the minor nature of the required modifications, we omit the proof of Lemma~\ref{Lemma_Car} here. For a more general version of the Carleman estimate, applicable when the Laplacian is replaced by a general elliptic operator, we refer the reader to \cite{NguyenLiKlibanov:2019}.

We have the theorem.	
	\begin{Theorem}[Carleman Convexification Theorem]
Let $\lambda_0$ be as given in Lemma \ref{Lemma_Car}. The following statements hold:
\begin{enumerate}
    \item For all $\lambda > 1$ and $\epsilon > 0$, the functional $J_{\lambda, \epsilon}$ is Fr\'echet differentiable, and its derivative $DJ_{\lambda, \epsilon}$ is Lipschitz continuous. That is, there exists a constant $L$, depending only on $\Omega$, $M$, and $N$, such that
    \[
        \|DJ_{\lambda, \epsilon}(\bv_2) - DJ_{\lambda, \epsilon}(\bv_1) \|_{{H^s(\Omega)^N}'} \leq L \|\bv_2 - \bv_1\|_{H^s(\Omega)^N}
    \]
    for all $\bv_1, \bv_2 \in B(0, M)$.

    \item There exists a constant $\lambda_1 = \lambda_1(M, N, r, \Omega) \geq \lambda_0$ such that for all $\epsilon > 0$ and $\lambda \geq \lambda_1$, the functional $J_{\lambda, \epsilon}$ is strictly convex in $B(0, M)$. Specifically,
    \begin{multline}
        J_{\lambda, \epsilon}(\bu) - J_{\lambda, \epsilon}(\bv) - DJ_{\lambda, \epsilon}(\bv)(\bu - \bv) 
        \geq (C_1\lambda - C_2)\int_{\Omega} e^{2\lambda (z - r)^2} \left( |\bu - \bv|^2 + |\nabla (\bu -  \bv)|^2 \right)\,d\x
        \\
        +
         (C_3\lambda^4 - C_4\lambda^3)\int_{\Gamma} e^{2\lambda (z - r)^2}  (|\bu - \bv|^2 +  |\partial_z (\bu - \bv)|^2)\,d\sigma(\x) 
                 + \epsilon \|\bu - \bv\|_{H^s(\Omega)^N}^2, \label{conv}
    \end{multline}
    for all $\bu, \bv \in H$, where $C_1, C_2, C_3$ and $C_4$ are positive constants depending only on $M$, $N$, $r$, and $\Omega$. As a result, $J_{\lambda, \epsilon}$ has a unique minimizer in $B(0, M)$, denoted by $\bv_{\rm min}$.

    \item Let $\lambda > \lambda_1$ and define $\Lambda = C_1\lambda - C_2 > 0$. Fix an initial guess $\bv^{(0)} \in B(0, M)$ and assume
    \[
        \left\{ \bm{\varphi} \in B(0, M) : \|\bm{\varphi} - \bv_{\rm min}\|_{H^s(\Omega)^N} < \|\bv^{(0)} - \bv_{\rm min}\|_{H^s(\Omega)^N} \right\}
        \subset B(0, M).    \]
    Define $\eta_0 = \min\{2\Lambda / L^2,\, 1\}$ and fix a step size $\eta \in (0, \eta_0)$. For each $m \geq 0$, set the iteration
    \begin{equation}
        \bv^{(m+1)} = \bv^{(m)} - \eta J_{\lambda, \epsilon}'(\bv^{(m)})
        \label{4.7}
    \end{equation}
    where $J_{\lambda, \epsilon}': H^s(\Omega)^N \to H^s(\Omega)^N$ is the Rietz representation of $DJ_{\lambda, \epsilon}$. That means,  
    	\[
		\langle J_{\lambda, \epsilon}'(\bv), \bm{\varphi}\rangle_{H^s(\Omega)^N} = DJ_{\lambda, \epsilon}(\bv)(\bm{\varphi})
		\quad 
		\mbox{for all } \bv, \bm{\varphi} \in H^s(\Omega)^N.
	\] 
    Then, there exists a constant $q \in (0, 1)$ such that for all $m \geq 0$,
    \[
        \bv^{(m)} \in B(0, M) \quad \text{and} \quad \|\bv^{(m)} - \bv_{\rm min}\|_{H^s(\Omega)^N} \leq q^{m-1}\|\bv^{(0)} - \bv_{\rm min}\|_{H^s(\Omega)^N}.
    \]
    Consequently, the sequence $\{\bv^{(m)}\}_{m \geq 0}$ converges to the unique minimizer $\bv_{\rm min}$ as $m \to \infty$.
\end{enumerate}
\label{Thm_conv}
\end{Theorem}

The main ideas to establish convexification theorems were originally introduced in \cite{KlibanovIoussoupova:SMA1995} and later applied to inverse scattering problems in \cite{KhoaKlibanovLoc:SIAMImaging2020, LeNguyen:JSC2022}. In those works, some versions of the convexification theorem were formulated for corresponding versions of the cost functional \( J_{\lambda, \epsilon} \) that does not include an integral term over the measurement surface \( \Gamma \) and the complementary surfaces $\partial \Omega \setminus \Gamma$. In the absence of these terms, the cost functional is defined on the set
\begin{multline*}
	\Big\{
		\bm{\varphi}
		= 
		\begin{bmatrix}
			\varphi_1 & \dots & \varphi_N
		\end{bmatrix}^{\rm T}
		\in H^s(\Omega)^N : 
		\|\bm{\varphi}\|_{H^s(\Omega)^N} < M, \ 
		\varphi_m|_\Gamma = g_m, \ 
		\partial_z \varphi_m|_{\Gamma} = h_m,
		\\ 
		\varphi_m|_{\partial \Omega \setminus \Gamma} = 0,\ m = 1, \dots, N
	\Big\}.
\end{multline*}
However, it is nontrivial to verify whether this set is nonempty, whereas the nonemptiness of the set \( H \) defined earlier is straightforward. To circumvent this difficulty, we incorporate the integral over \( \Gamma \) and $\partial \Omega \setminus \Gamma$ into the cost functional \( J_{\lambda, \epsilon} \).

The proof of Theorem~\ref{Thm_conv} closely follows the methodology presented in \cite[Theorem 4.1]{LeLeNguyen:2024}. The first part is derived through straightforward algebraic manipulations. The second part builds upon earlier results from \cite{KlibanovNik:ra2017}, \cite[Theorem 5.1]{KhoaKlibanovLoc:SIAMImaging2020}, and \cite[Theorem 4.1, part 2]{LeLeNguyen:2024}, with a minor but important modification: in contrast to these prior works, where expressions such as \( (C_1\lambda - C_2) \) and \( (C_3\lambda^4 - C_4\lambda^3) \) are replaced with a generic constant \( C \), we preserve the explicit dependence on \( \lambda \) to emphasize the crucial role played by the Carleman parameter in the convexity inequality \eqref{4.2}. This modification is justified by adapting the proof technique in \cite[Theorem 4.1]{LeLeNguyen:2024}, with the key difference being the use of the Carleman estimate in \eqref{4.55} in place of the one used in that reference.
 The third part of Theorem~\ref{Thm_conv} follows directly from \cite[Theorem 2]{LeNguyen:JSC2022}.

We next discuss how close the minimizer of $J_{\lambda, \epsilon}$ is to the true solution of \eqref{4.2}. 
Let \[ \bv^* = \left[
	\begin{array}{ccc}
		v^*_1 & \dots & v^*_N
	\end{array}
\right]^{\rm T} \] denote the exact solution to \eqref{4.2} corresponding to the noiseless boundary data \( g_m^* \) and \( h_m^* \), which are idealized versions of the measured data \( g_m \) and \( h_m \), respectively. That is, for each \( m \in \{1, \dots, N\} \), the functions \( v_m^* \) satisfy the following boundary value problem:
\begin{equation}
    \left\{
        \begin{array}{ll}
            \Delta v_m^*(\x) 
            + \sum\limits_{n = 1}^N \nabla v_n^*(\x) \cdot \widetilde{\bf b}_{mn}		 
            + \sum\limits_{n = 1}^N \sum\limits_{l = 1}^N \widetilde{a}_{mnl} \nabla v_n^*(\x) \cdot \nabla v_l^*(\x) 
            = 0, & \x \in \Omega, \\
            v_m^*(\x) = g_m^*(\x), & \x \in \Gamma, \\
            \partial_z v_m^*(\x) = h_m^*(\x), & \x \in \Gamma, \\
            v_m^*(\x) = 0, & \x \in \partial \Omega \setminus \Gamma.
        \end{array}
    \right.
    \label{4.5}
\end{equation}
Additionally, we assume the noise level in the measured data is bounded by $\delta$ in the following sense
\begin{equation}
\sum_{m = 1}^N \big(\|g_m - g^*\|_{H^1(\Gamma)} + \|h_m - h_m^*\|_{L^2(\Omega)}\big) < \delta,
\end{equation}
for some small constant \( \delta \ll 1 \). Then, using the same arguments of \cite[Theorem 4.2]{LeLeNguyen:2024}, we have
\begin{equation}
	\|\bv_{\rm min} - \bv^*\|_{H^1(\Omega)^N} \leq C\big(\delta + \sqrt{\epsilon}\|\bv^*\|_{H^s(\Omega)^N}\big)
	\label{4.9}
\end{equation}
for some constant $C$ depending only on $\Omega, M$, $r$, $\lambda$, and $N$.
In \eqref{4.9}, we employ the standard Sobolev norm instead of incorporating a Carleman-weighted functional, as was done in our previous work on convexification. This formulation remains valid under the assumption that \( \lambda \) is fixed and the constant \( C \) may depend on \( \lambda \). A direct implication of \eqref{4.9} is that the global minimizer of \( J_{\lambda, \epsilon} \) provides an approximation to the true solution of \eqref{4.2}, with an error bounded by \( O(\delta + \sqrt{\epsilon}) \).

The WKB method for phase retrieval presented in Section~\ref{sec_phase}, the derivation of the frequency dimensional reduction model in Section~\ref{sec_reduceEqn}, the convexification result in Theorem~\ref{Thm_conv}, and the error estimate in \eqref{4.9} collectively motivate the design of Algorithm~\ref{alg} for solving the phaseless inverse scattering problem.

\section{Numerical study} \label{sec_num}
In this section, we highlight key aspects of Algorithm~\ref{alg}'s implementation and present several numerical examples to demonstrate its performance.

\begin{algorithm}[h!]
\caption{\label{alg}The Carleman convexification method to compute the numerical solution to the phaseless inverse scattering problem}
	\begin{algorithmic}[1]
	\State   Having the data $f$ in hand, for each $k$, we minimize the mismatch functional $J_k$ defined in \eqref{2.5} using the initial guess $u_{\rm init}(\x, k) = f(\x, k) e^{\ik |\x - \x_0|}$ as in \eqref{2.5555}. The obtained minimizer $u(\x, k)$, $(\x, k) \in \Gamma_L \times [\underline k, \overline k]$ is the desired wave function including the phase information. \label{s1}
	\State Compute $u(\x, k)$ and $\partial_z u(\x, k)$ on $\Gamma,$ which is the bottom portion of $\partial \Omega.$
	\State  Choose a cut-off number $N$, a Carleman parameter $\lambda$, and a regularization parameter $\epsilon$.
	 Choose Carleman parameters $\x_0,$ $\beta$, and $\lambda$ and a regularization parameter $\epsilon$.
	 \State Minimize the strictly convex Carleman weighted functional $J_{\lambda, \epsilon}$ defined in \eqref{4.4}. The minimizer is denoted by $\bv_{\rm comp}(\x) = \left[
	 	\begin{array}{cccc}
		v_1^{\rm comp} & v_2^{\rm comp} &\dots & v_N^{\rm comp}
		\end{array}
	 \right]^{\rm T}$, $\x \in \Omega.$ \label{s4}
	\State By \eqref{3.6}, we compute $v^{\rm comp}(\x, k)$ using \begin{equation}
	v^{\rm comp}(\x, k) = \sum_{n = 1}^N v_n^{\rm comp}(\x) \Psi_n(k),
\end{equation}
for $\x \in \Omega,$ $k \in [\underline k, \overline k]$.
	\State 
	Due to \eqref{3.4}, a numerical solution to Problem \ref{cip} can be computed via
	\begin{multline*}
		c^{\rm comp}(\x) = 1 - \frac{1}{\overline k - \underline k}\int_{\underline k}^{\overline k}\Re e\Bigg[\Delta v^{\rm comp}(\x, k) + k^2 [\nabla v^{\rm comp}(\x, k)]^2 
		\\
		+ 2 \Big(\ik - \frac{1}{|\x - \x_0|}\Big) \nabla v^{\rm comp}(\x, k) \cdot \frac{\x - \x_0}{|\x - \x_0|}\Bigg]d\x d\theta  
	\end{multline*}
	for all $\x \in \Omega.$	 
\end{algorithmic}
\end{algorithm}

To generate the simulated data, we set \( R = 1 \), placing the source at \( \x_0 = (0, 0, -4) \) and using the wave number interval \( [\pi, 2\pi] \). To generate simulated data, the domain \( \Omega = (-1, 1)^3 \) is discretized using a uniform grid defined as
\begin{equation*}
	\mathcal G = \left\{
		(x_i = -1 + (i - 1)d_\x,\, y_j = -1 + (j - 1)d_\x,\, z_t = -1 + (t - 1)d_\x): 1 \leq i, j, t \leq N_\x
	\right\},
\end{equation*}
where \( N_\x = 21 \) and \( d_\x = 2/(N_\x - 1) \). The wave number is  interval set to be $[\pi, 2\pi]$ and is  uniformly discretized into
\[
	\mathcal K = \{k_1 = \underline k,\, k_2,\, \dots,\, k_{N_k} = \overline k\},
\]
with \( k_i = \underline k + (i - 1)\frac{\overline k - \underline k}{N_k - 1} \) and \( N_k = 121 \). The forward problem is addressed by reformulating the Helmholtz model \eqref{Hel} into the Lippmann-Schwinger integral equation, following the approach in \cite{ColtonKress:2013}. This integral equation is then numerically solved using the volume integral equation method developed in \cite{LechleiterNguyen:acm2014, Nguyen:anm2015}. Let \( u^*(\x, k) \), with \( \x \in \mathcal G \) and \( k \in \mathcal K \), denote the exact solution. 
We define the noisy data as \[
	f(\x, k) = |u^*(\x, k)|(1 + \delta\, \mathrm{rand}),
\]
for \( (\x, k) \in (\Gamma_L \cap \mathcal G) \times \mathcal K \), where \( \delta = 10\% \) and \( \mathrm{rand} \) represents a uniformly distributed random variable in the interval \( [-1, 1] \).  
In this section, we set $L = 0.28$.
For the artificial parameters involved in solving the inverse problem, we select \( N = 7 \), \( \lambda = 1.1 \), and \( \epsilon = 10^{-5.75} \). These values are determined through a trial-and-error process. They are consistent with the corresponding set of parameters in \cite{LeNguyen:JSC2022}.  Specifically, we use Test 1 as a reference case to manually identify a suitable set of parameters, which are then consistently applied across all subsequent tests. 

\subsection{Numerical implementation of phase retrieval}

\begin{figure}[h!]
    \flushright
    \subfloat[\label{fig2a}Noisy data $f$]{\includegraphics[width=0.31\textwidth]{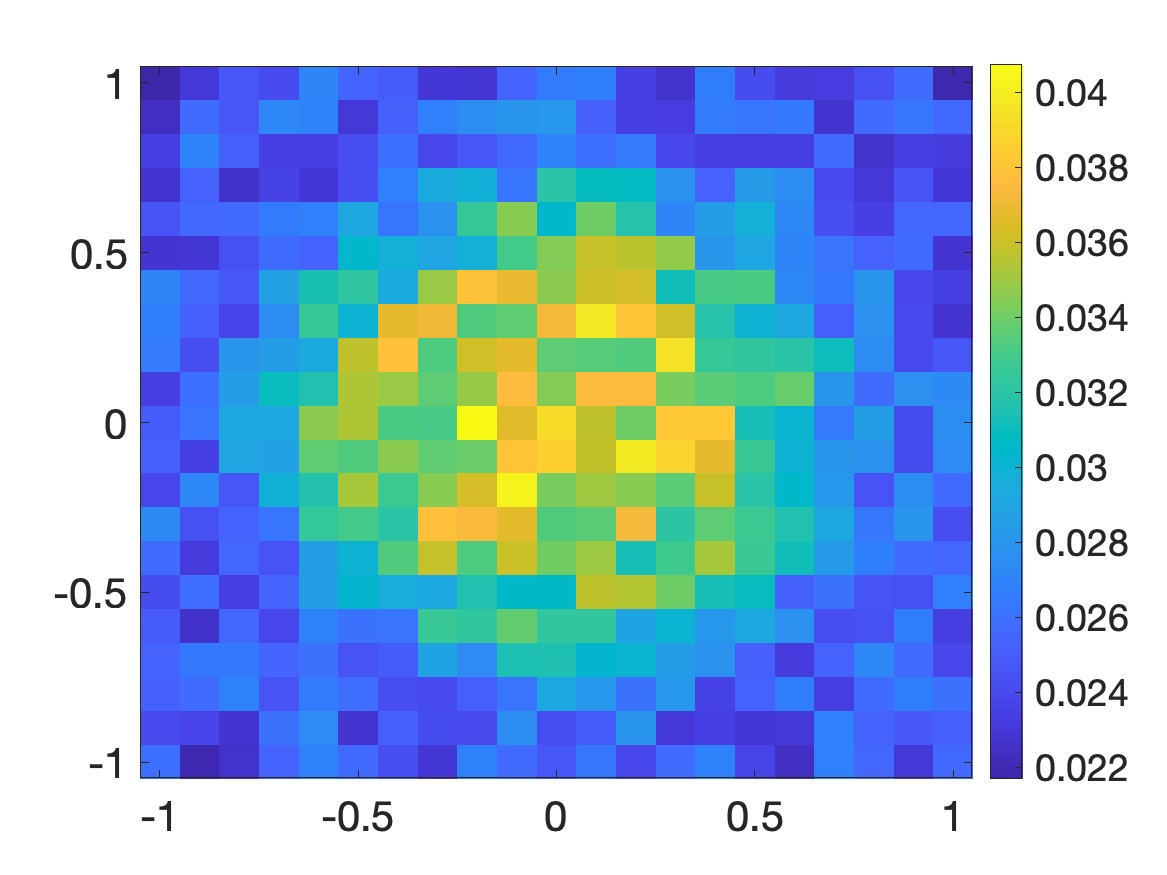}}
    \quad
    \subfloat[\label{fig2b}True function $\Re e(u^*)$]{\includegraphics[width=0.31\textwidth]{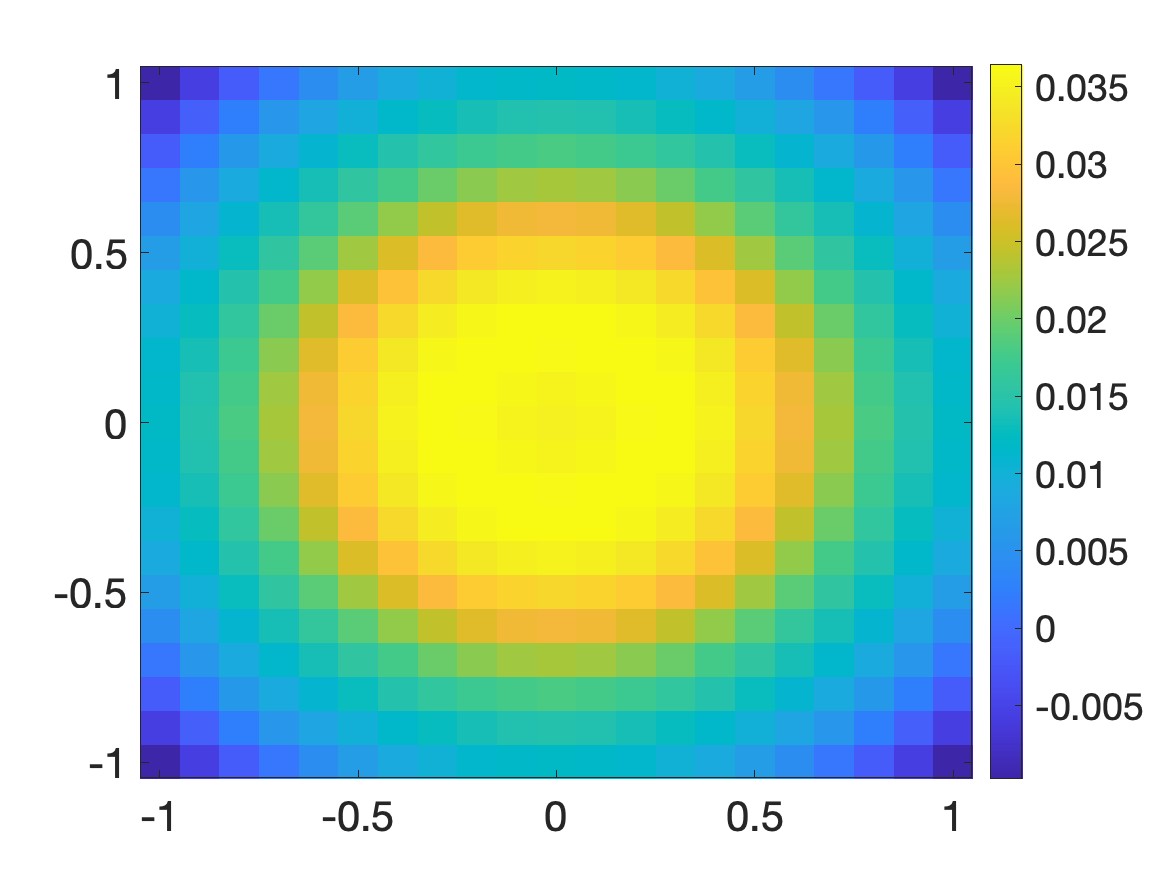}}
    \quad
    \subfloat[\label{fig2c}Reconstructed function $\Re e(u)$]{\includegraphics[width=0.31\textwidth]{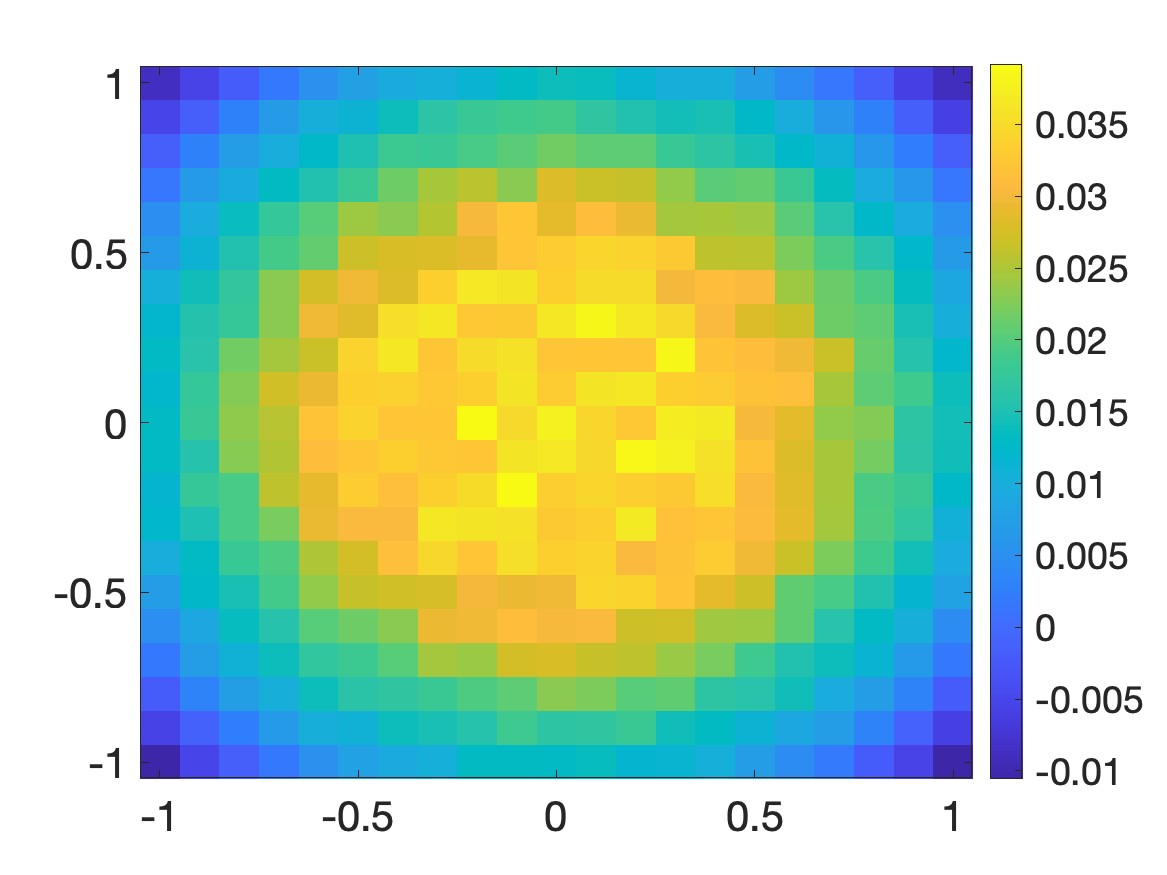}}
    
    \subfloat[\label{fig2d}True function $\Im m(u^*)$]{\includegraphics[width=0.32\textwidth]{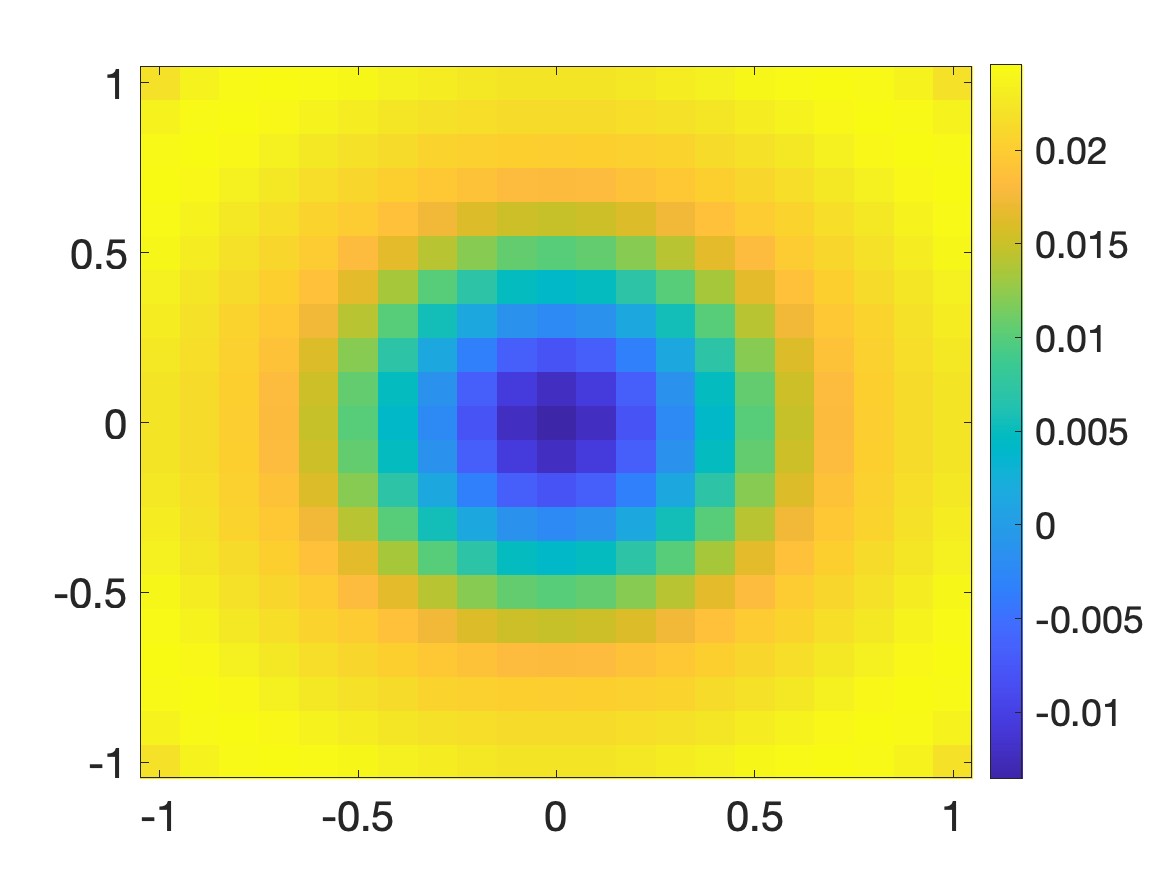}}
    \quad
    \subfloat[\label{fig2e}Reconstructed function $\Im m(u)$]{\includegraphics[width=0.31\textwidth]{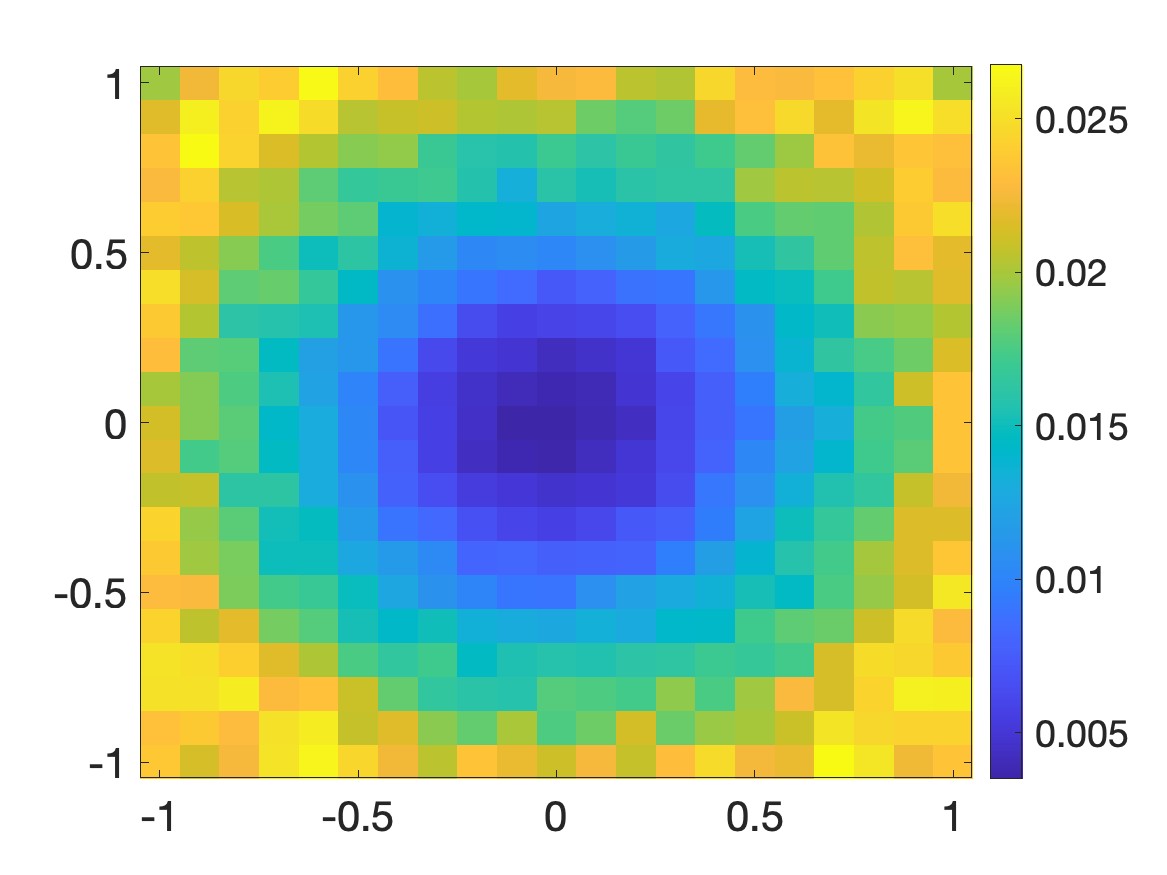}}
    \caption{\label{fig_phase} Reconstruction of the wave function \( u(\x, k) \) on \( \Gamma \) with \( k = \overline{k} \), showing both real and imaginary parts from magnitude-only data. The data in (a) corresponds to Test 1 with a noise level of 10\%.}
\end{figure}

In Step~\ref{s1} of Algorithm~\ref{alg}, the cost functional \( J_k \) is approximated using a Riemann sum over the spatial grid \( \Gamma_L \cap \mathcal{G} \) and the frequency partition \( \mathcal{K} \). To minimize \( J_k \), we use MATLAB's built-in optimization routine \texttt{fminunc}, which proves to be both efficient and effective in practice. As a representative example, we demonstrate the recovery of the lost phase. Starting from the noisy data $f$ in Test 1, the optimization is initialized with \( u_{\text{init}} = f(\x, k) e^{i k |\x - \x_0|} \), and \texttt{fminunc} is applied to minimize the cost. The resulting reconstruction of the complex-valued wave function \( u \), including both its real and imaginary parts, is presented in Figure~\ref{fig_phase}.

The numerical results in Figure~\ref{fig_phase} demonstrate the effectiveness of the phase retrieval strategy. Despite the presence of $10\%$ noise in the measured data (Figure~\ref{fig2a}), the reconstructed wave field (Figures~\ref{fig2c} and~\ref{fig2e}) closely approximates the true solution (Figures~\ref{fig2b} and~\ref{fig2d}) in both the real and imaginary components. The method successfully preserves key structural features, such as spatial distribution and amplitude. In particular, the real part of the reconstruction captures the peak profile with good accuracy, while the imaginary part accurately recovers the expected dip pattern.

Although the input data in Figure~\ref{fig2a} appears heavily degraded, the reconstructed real and imaginary parts (Figures~\ref{fig2c} and~\ref{fig2e}) exhibit a clear improvement in quality. This enhancement is primarily attributed to the structure of the cost functional \( J_k \), which incorporates a regularization term involving the Laplacian of the wave field (see the first term on the right-hand side of \eqref{2.5}). This term promotes smoothness by penalizing irregular or non-differentiable behavior, resulting in reconstructions that are inherently smoother than the raw, noisy magnitude data.
Despite this regularization, some residual noise persists in the reconstructed fields. This is due to the high noise level (10\%) in the input and the effect of the data fidelity term in \eqref{2.5}, which constrains the solution to remain close to the measured data \( f \), thereby preserving some of its noise features. To address this, a Fourier filtering procedure (see \eqref{3.6}) is employed in a subsequent step to suppress remaining high-frequency components, further improving the stability and accuracy of the final reconstruction.

\subsection{Numerical implementation of the convexification method}

The Carleman convexification framework is implemented in Step~\ref{s4} of Algorithm~\ref{alg}. The goal of this step is to compute the vector \( \bv_{\rm comp} \), whose components represent the Fourier coefficients of the function \( v \) defined in \eqref{change_variable}. This vector corresponds to the solution of problem \eqref{3.11}, which is equivalent to \eqref{4.2}. However, it is important to note that \eqref{4.2} involves the inverse matrix \( S^{-1} \), and numerical observations indicate that some of its entries can attain large magnitudes, potentially degrading the accuracy of the final reconstruction. To mitigate this issue and avoid unnecessary numerical errors, we opt to solve \eqref{3.11} directly rather than \eqref{4.2}. 

In this case, we reformulate the cost functional \( J_{\lambda, \epsilon} \) as follows:
\begin{multline}
	J_{\lambda, \epsilon}(\bm{\varphi}) = \sum_{m = 1}^N \Bigg[\int_{\Omega} e^{2\lambda (z - r)^2} \Big|
		\sum_{n = 1}^N s_{mn} \Delta \varphi_n(\x) 
	 +  \sum_{n = 1}^N\sum_{l = 1}^N a_{mnl}\nabla \varphi_n(\x) \cdot \nabla \varphi_l(\x) 	
	 \\ 	
	+   \sum_{n = 1}^N  {\bf b}_{mn} \cdot \nabla \varphi_n(\x)
	\Big|^2 d\x
	+ \lambda^4 e^{2\lambda (R+r)^2}\int_{\Gamma} \big( |\varphi_m - g_m|^2 + |\partial_z \varphi_m - h_m|^2\big)d\sigma(\x)
	\\
		+ \lambda^4 \int_{\partial \Omega \setminus \Gamma} e^{2\lambda(z- r)^2} |\varphi_m|^2 d\sigma(\x)
	+ \epsilon \|\varphi_m\|_{H^2(\Omega)}^2\Bigg]
	\label{5.2}
\end{multline}
for all \( \bm{\varphi} \in H \). Without confusion, we continue to denote this functional by \( J_{\lambda, \epsilon} \).
We also note that, unlike in the theoretical part, which requests $s \geq 4$, we use the regularization norm $H^2(\Omega)$. This significantly simplifies the implementation without reducing the quality of the final reconstruction.

Minimizing the functional \( J_{\lambda, \epsilon} \) in \eqref{5.2} requires an initial guess \( \bv^{(0)} \). According to the convexification theorem, this initial guess does not need to be close to the global minimizer of \( J_{\lambda, \epsilon} \); the only requirement is that \( \bv^{(0)} \in H \). Following the approach in \cite{LeNguyen:JSC2022}, we simplify the construction of \( \bv^{(0)} \) by omitting the nonlinear term \( \sum_{n = 1}^N\sum_{l = 1}^N a_{mnl}\nabla \varphi_n(\x) \cdot \nabla \varphi_l(\x) \) from \eqref{5.2}. That means, we  define \( \bv^{(0)} \) as the minimizer of the simplified strictly convex cost functional:
\begin{multline}
	J^{(0)}_{\lambda, \epsilon}(\bm{\varphi}) = \sum_{m = 1}^N \Bigg[\int_{\Omega} e^{2\lambda (z - r)^2} \Big|
		\sum_{n = 1}^N s_{mn} \Delta \varphi_n(\x) 
	+   \sum_{n = 1}^N  {\bf b}_{mn} \cdot \nabla \varphi_n(\x)
	\Big|^2 d\x
	+ \lambda^4 \int_{\partial \Omega \setminus \Gamma} e^{2\lambda(z- r)^2} |\varphi_m|^2 d\sigma(\x)
	\\
	+ \lambda^4 e^{2\lambda (R+r)^2}\int_{\Gamma} \big( |\varphi_m - g_m|^2 + |\partial_z \varphi_m - h_m|^2\big)d\sigma(\x)
	+ \epsilon \|\varphi_m\|_{H^2(\Omega)}^2\Bigg]
	\label{5.3}
\end{multline}
for all \( \bm{\varphi} \in B(0, M) \). 
Having the initialized vector $\bv^{(0)}$ in hand, we compute the minimizer by the gradient descent method using the formula \eqref{4.7}. 
This requires us to compute the derivative of $J_{\lambda, \epsilon}$. In the implementation, we compute the derivative in finite difference by regarding $J_{\lambda, \epsilon}$ as a function of $N_\x^3 N$ variables. 
The discretized version of $\bm{\varphi} \in H^2(\Omega)^N$ is 
\[
	\bm{\varphi} := \big\{\varphi_m(x_i, y_j, z_t): 1 \leq i, j, t \leq N_x, 1 \leq m \leq N\big\}
\]
and that of $J_{\lambda, \epsilon}$ is 
\begin{multline}
	J_{\lambda, \epsilon}(\bm{\varphi}) = \sum_{m = 1}^N \Bigg[d_\x^3\sum_{i, j, t = 2}^{N_\x-1} e^{2\lambda (z_t - r)^2} \Big|
		\sum_{n = 1}^N s_{mn} \Delta^{d_\x} \varphi_n(x_i, y_j, z_t) 
		\\
	 +  \sum_{n = 1}^N\sum_{l = 1}^N a_{mnl}\nabla^{d_\x} \varphi_n(x_i, y_j, z_t) \cdot \nabla^{d_\x} \varphi_l(x_i, y_j, z_t)
	+   \sum_{n = 1}^N  {\bf b}_{mn} \cdot \nabla^{d_\x} \varphi_n(x_i, y_j, z_t)
	\Big|^2 
	\\
	+ \lambda^4 e^{2\lambda (R+r)^2} d_\x^2\sum_{i, j = 1}^{N_\x} \big( |\varphi_m(x_i, y_j, z_1) - g_m(x_i, y_j, z_1)|^2 + |\partial_z^{d_\x} \varphi_m(x_i, y_j, z_1) - h_m(x_i, y_j, z_1)|^2\big)
	\\
		+ \lambda^4 d_{\x}^2 \sum_{i, j = 1}^{N_\x} e^{2\lambda(z_{N_\x}- r)^2} |\varphi_m(x_i, y_j, z_{N_\x})|^2 d\sigma(\x)
		+  \lambda^4 d_{\x}^2 \sum_{i, j \in\{1, N_\x\}} \sum_{t = 1}^{N_\x} e^{2\lambda(z_t- r)^2} |\varphi_m(x_i, y_j, z_t)|^2 d\sigma(\x)
		\\
	+ \epsilon d_\x^3 \sum_{i, j, t = 2}^{N_\x-1} |\varphi_m(x_i, y_j, z_t)|^2 + |\nabla^{d_\x}\varphi_m(x_i, y_j, z_t)|^2 + |\Delta^{d_\x}\varphi_m(x_i, y_j, z_t)|^2\Bigg].
	\label{5.4}
\end{multline}

In \eqref{5.4}, 
\begin{align*}
	\nabla^{d_\x} \varphi(x_i, y_j, z_t) 
	&=\left[
		\begin{array}{c}
			\partial_x^{d_\x}\varphi(x_i, y_j, z_t)\\
			\partial_y^{d_\x}\varphi(x_i, y_j, z_t)
			\\
			\partial_z^{d_\x}\varphi(x_i, y_j, z_t)
		\end{array}
	\right]
	= 
	\left[
		\begin{array}{c}
			\frac{\varphi(x_{i+1}, y_j, z_t) - \varphi(x_{i-1}, y_j, z_t)}{2d_\x}\\
			\frac{\varphi(x_{i}, y_{j+1}, z_t) - \varphi(x_{i}, y_{j-1}, z_t)}{2d_\x}\\
			\frac{\varphi(x_{i}, y_j, z_{t+1}) - \varphi(x_{i}, y_j, z_{t-1})}{2d_\x}
		\end{array}
	\right],
	\\
	\Delta^{d_\x}   \varphi(x_i, y_j, z_t)
	&= \frac{1}{d_\x^2}\Big[\varphi(x_{i+1}, y_j, z_t) + \varphi(x_{i-1}, y_j, z_t) + \varphi(x_i, y_{j+1}, z_t)\\
	&\quad + \varphi(x_i, y_{j-1}, z_t) + \varphi(x_i, y_j, z_{t+1}) + \varphi(x_i, y_j, z_{t-1}) -6 \varphi(x_i, y_j, z_t)\Big].
\end{align*}
We can interpret the discretized form in \eqref{5.4} as a polynomial in the variables \[ \left\{\varphi_m(x_i, y_j, z_t) : 1 \leq i, j, t \leq N_\x,\ 1 \leq m \leq N \right\},\] allowing us to compute its derivative explicitly.

\subsection{Numerical examples}

We present three numerical solutions to the phaseless inverse scattering problem due to Algorithm \ref{alg}. 

\subsubsection{Test 1}

\begin{figure}[h!]
\centering
	\subfloat[\label{test1a}3D view of the true scatterer]{\includegraphics[width=.45\textwidth]{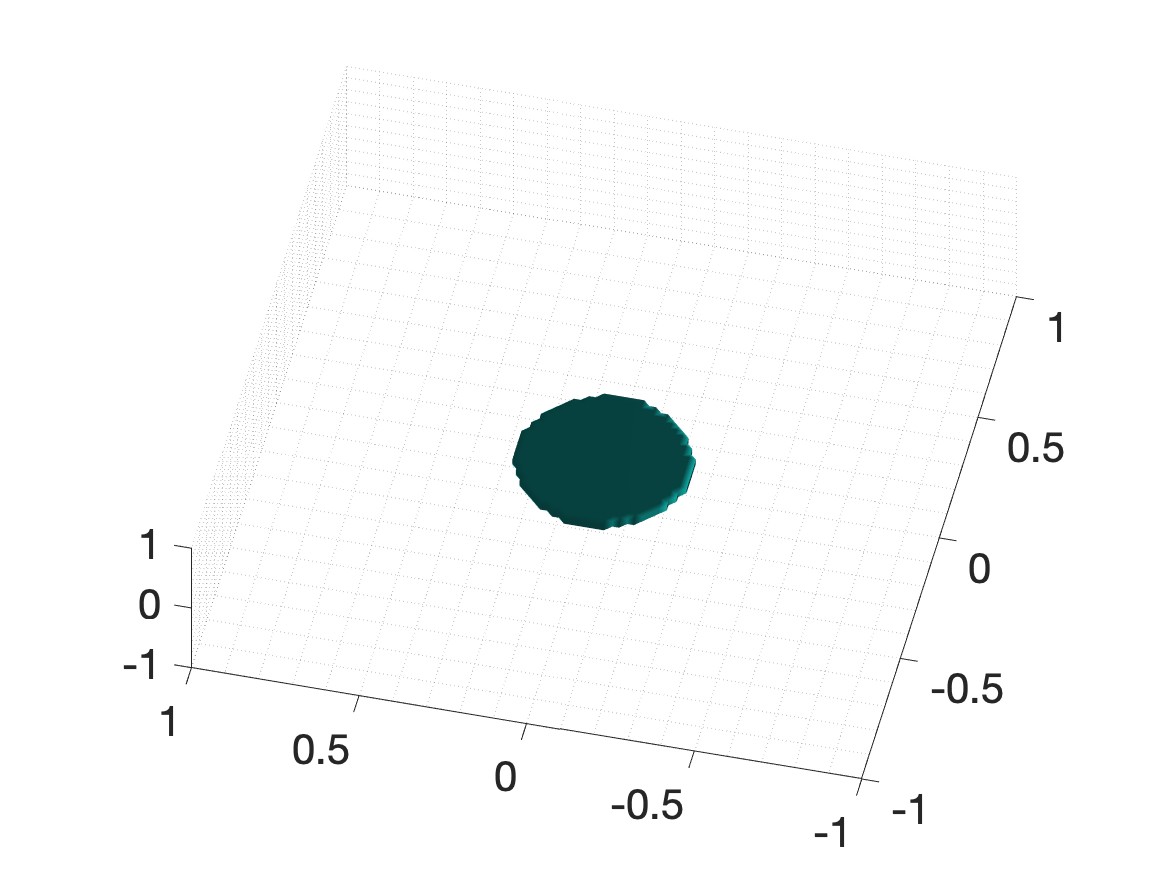}}
	\quad
	\subfloat[\label{test1b}3D view of the reconstructed scatterer]{\includegraphics[width=.45\textwidth]{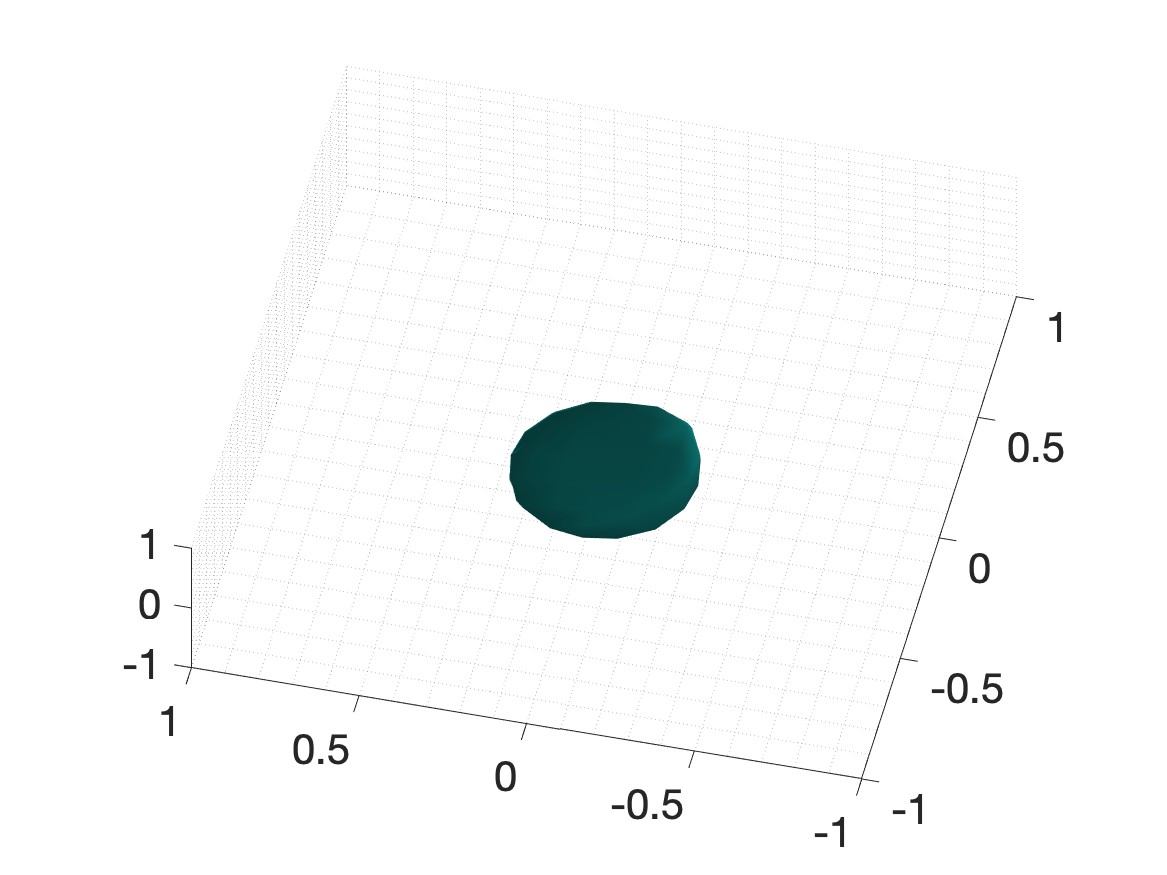}}

	\subfloat[\label{test1c}Cross-sectional view of the function \( c^{\rm true} \)]{\includegraphics[width=.45\textwidth]{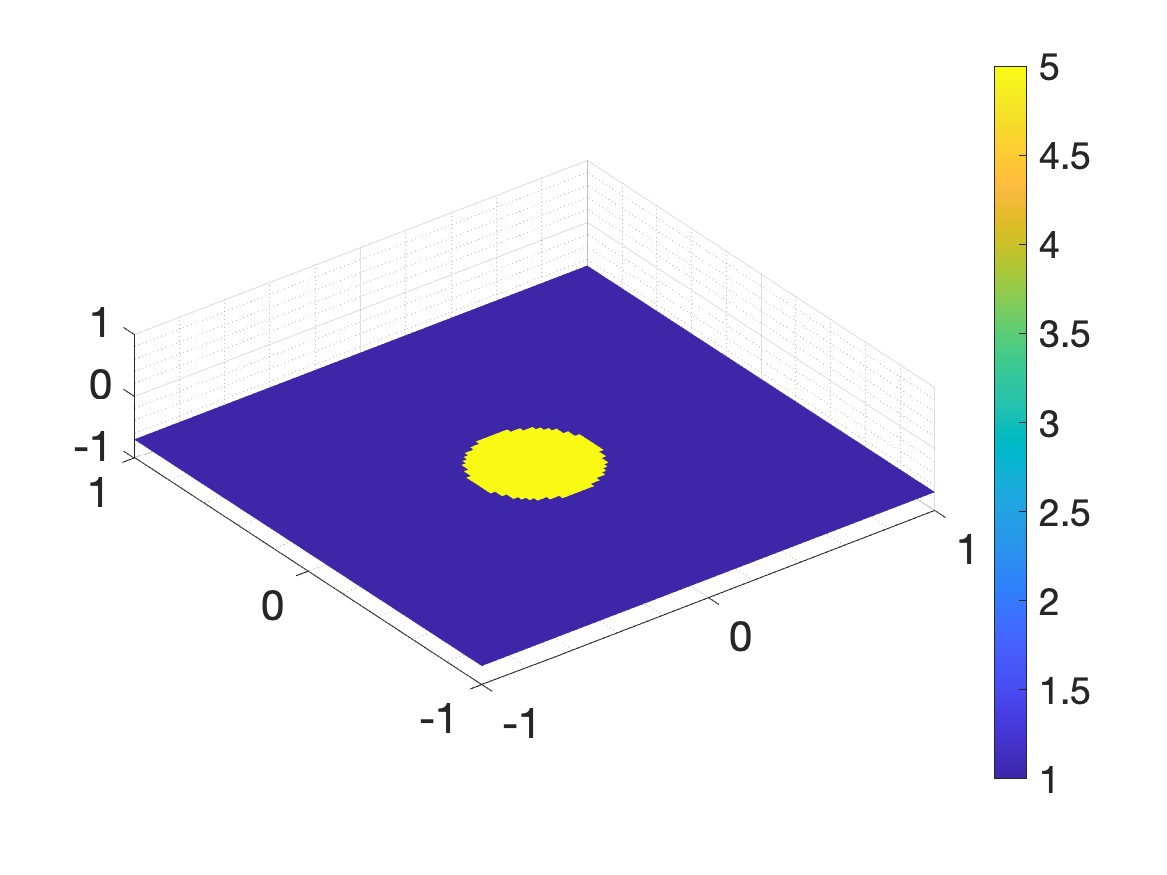}}	
	\quad
	\subfloat[\label{test1d}Cross-sectional view of the function \( c^{\rm comp} \)]{\includegraphics[width=.45\textwidth]{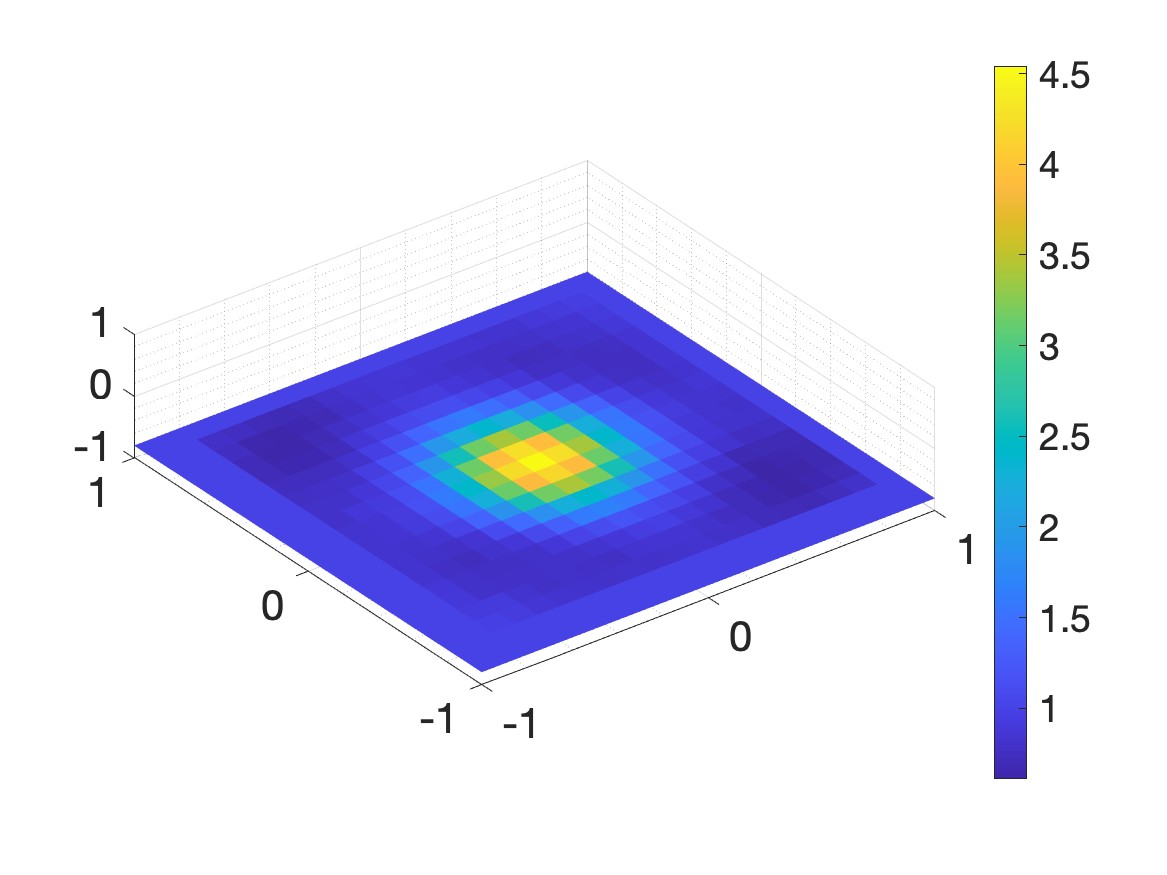}}

	\caption{\label{test1}Test 1. Visualization of the true and reconstructed dielectric constant \( c(\x) \), \( \x \in \Omega \). Subfigures (a) and (b) show 3D views of the true and reconstructed scatterers, respectively. Subfigures (c) and (d) display cross-sectional slices of the corresponding dielectric profiles. The reconstruction demonstrates strong agreement with the true target in both shape and amplitude. The dielectric function \( c(\x) \) was reconstructed from phaseless data containing 10\% noise.}
\end{figure}

We define the true profile of the dielectric constant as
\[
	c^{\rm true}(x, y, z) = \left\{
		\begin{array}{ll}
			5 & \text{if } x^2 + y^2 < 0.25^2 \text{ and } |z + 0.65| < 0.05,\\
			1 & \text{otherwise}.
		\end{array}
	\right.
\]
The true and reconstructed dielectric profiles are visualized in Figure~\ref{test1}.

Figure~\ref{test1} illustrates the effectiveness of the proposed reconstruction method in recovering the dielectric constant \( c(\x) \) from phaseless measurements corrupted with 10\% noise. The 3D visualizations in subfigures~\ref{test1a} and~\ref{test1b} show that the reconstructed scatterer closely matches the true inclusion in both position and geometry. Additionally, the cross-sectional views in subfigures~\ref{test1c} and~\ref{test1d} demonstrate that the reconstruction accurately preserves the spatial distribution and peak amplitude of the dielectric profile. The reconstructed function \( c \) achieves a maximum value of 4.54, corresponding to a relative error of 9.17\%. These results highlight the robustness and accuracy of the method, confirming its potential to produce high-quality reconstructions even from limited and noisy phaseless data.

\subsubsection{Test 2}

\begin{figure}[h!]
\centering
	\subfloat[\label{test2a}3D view of the true scatterer]{\includegraphics[width=.45\textwidth]{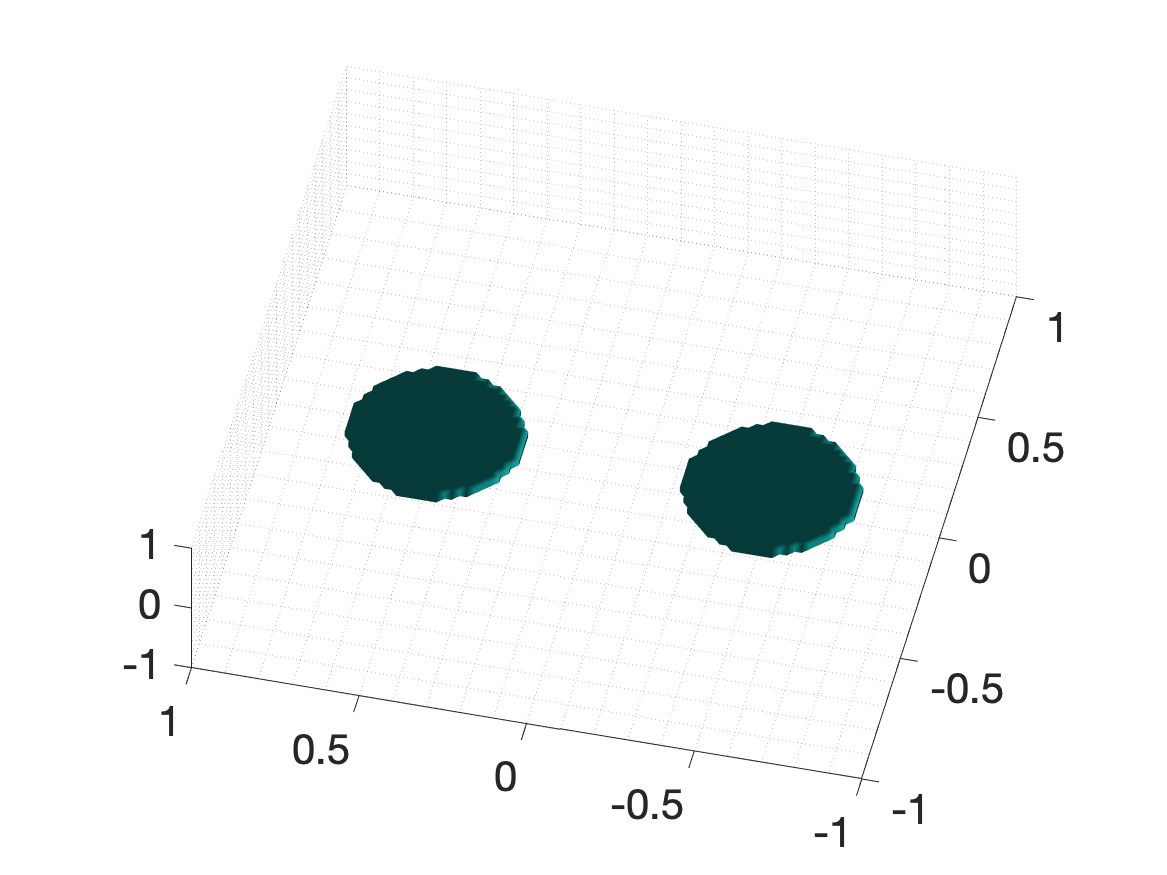}}
	\quad
	\subfloat[\label{test2b}3D view of the reconstructed scatterer]{\includegraphics[width=.45\textwidth]{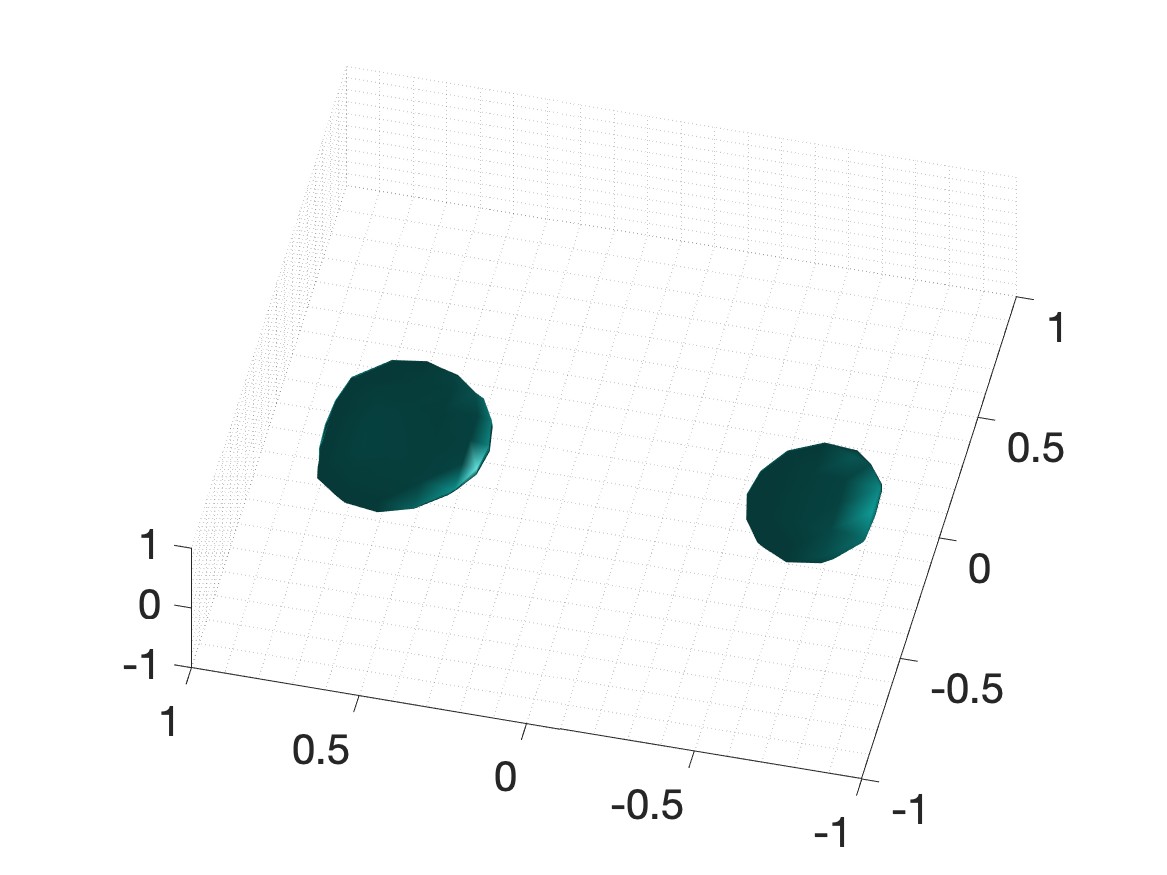}}

	\subfloat[\label{test2c}Cross-sectional view of the function \( c^{\rm true} \)]{\includegraphics[width=.45\textwidth]{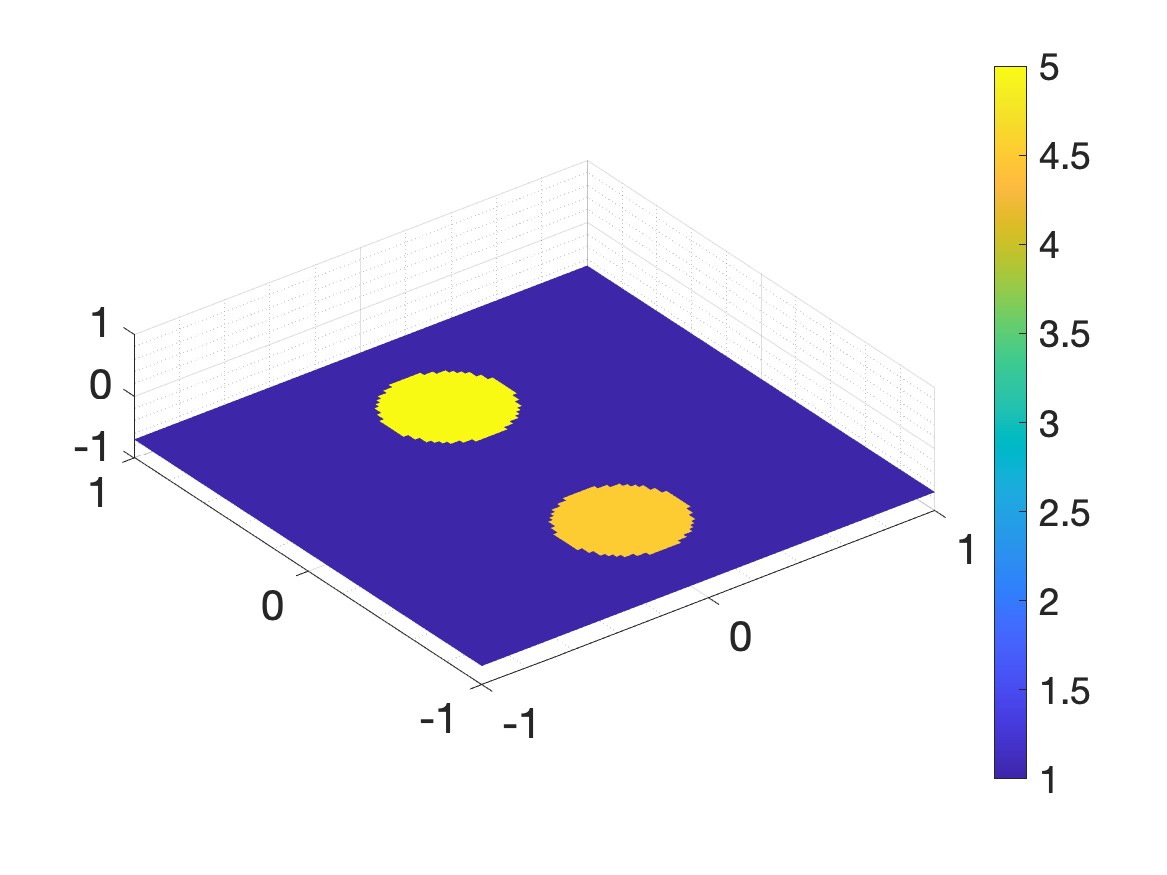}}	
	\quad
	\subfloat[\label{test2d}Cross-sectional view of the function \( c^{\rm comp} \)]{\includegraphics[width=.45\textwidth]{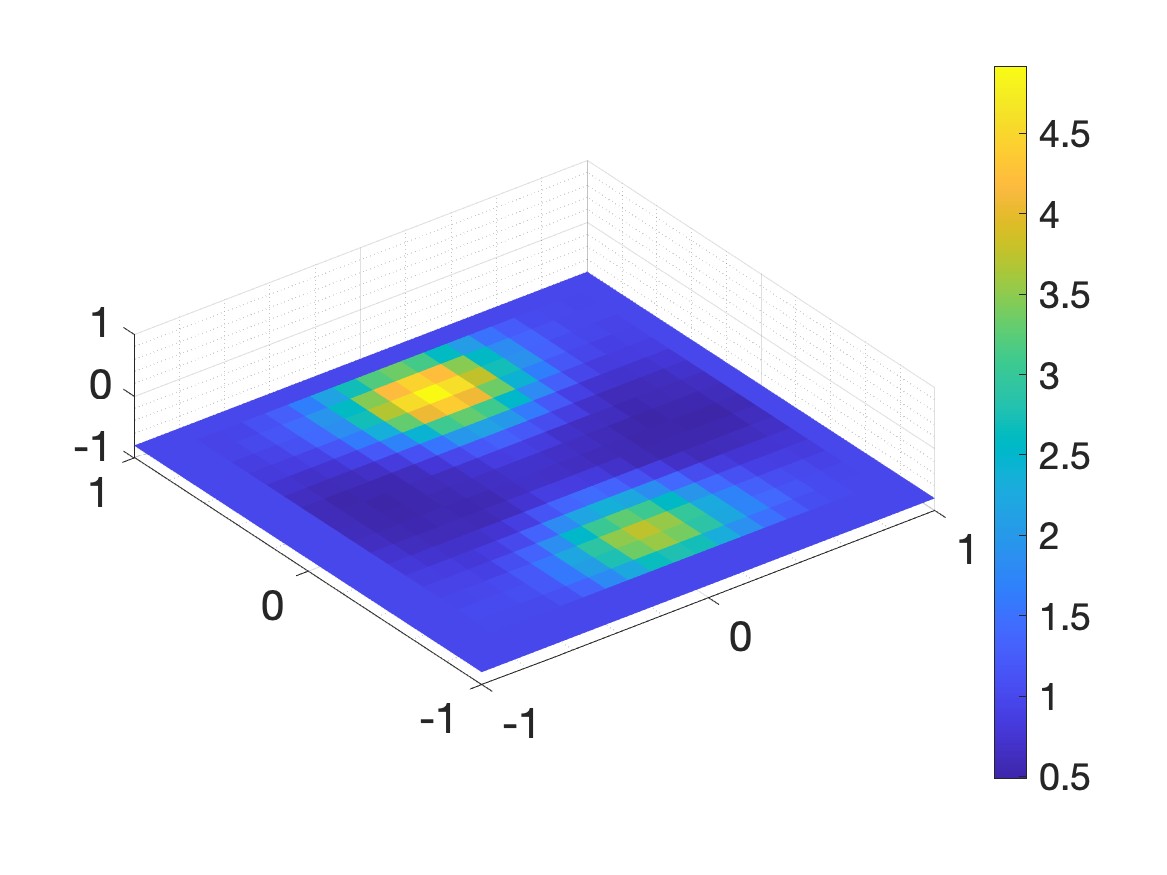}}

	\caption{\label{test2}Test 2. The true and reconstructed dielectric constant \( c(\x) \), \( \x \in \Omega \). Subfigures (a) and (b) display 3D views of the true and reconstructed scatterers, while (c) and (d) present cross-sectional views of their respective dielectric profiles. The reconstructed result successfully captures both the spatial location and amplitude of the true inclusions, demonstrating strong robustness to noise. The reconstruction was performed using phaseless data corrupted with 10\% noise.}
\end{figure}

We next test the case of two inclusions. The true dielectric constant function is given by
\[
	c^{\rm true} = 
	\left\{
		\begin{array}{ll}
			5 &(x - 0.5)^2 + y^2 < 0.25^2 \mbox{and} |z + 0.65| < 0.05,\\
			4.5 &(x + 0.5)^2 + y^2 < 0.25^2 \mbox{and} |z + 0.65| < 0.05,\\
			1 &\mbox{otherwise}.
		\end{array}
	\right.
\]
The true and reconstructed dielectric profiles are visualized in Figure~\ref{test2}.

Figure~\ref{test2} highlights the effectiveness of the proposed method in reconstructing multiple inclusions from phaseless data corrupted by 10\% noise. The 3D isosurface visualizations in subfigures~\ref{test2a} and~\ref{test2b} demonstrate that the reconstructed scatterers closely replicate the true geometry and spatial distribution of the targets. Importantly, the method accurately preserves the separation and relative sizes of the two inclusions. The cross-sectional views in subfigures~\ref{test2c} and~\ref{test2d} further validate the reconstruction quality, revealing that both the amplitude and localization of the high-contrast regions are well recovered. Quantitatively, the maximum value of the first inclusion centered at \( (0.5, 0, -0.65) \) is 4.92, yielding a relative error of 1.61\%, while the second inclusion centered at \( (-0.5, 0, -0.65) \) attains a maximum of 3.75, with a relative error of 16.78\%. These reconstruction errors are considered acceptable, especially given the high noise level and the severe ill-posedness of the inverse problem with single-sided measurements. 

\subsubsection{Test 3}

\begin{figure}[h!]
\centering
	\subfloat[\label{test3a}3D view of the true scatterer]{\includegraphics[width=.45\textwidth]{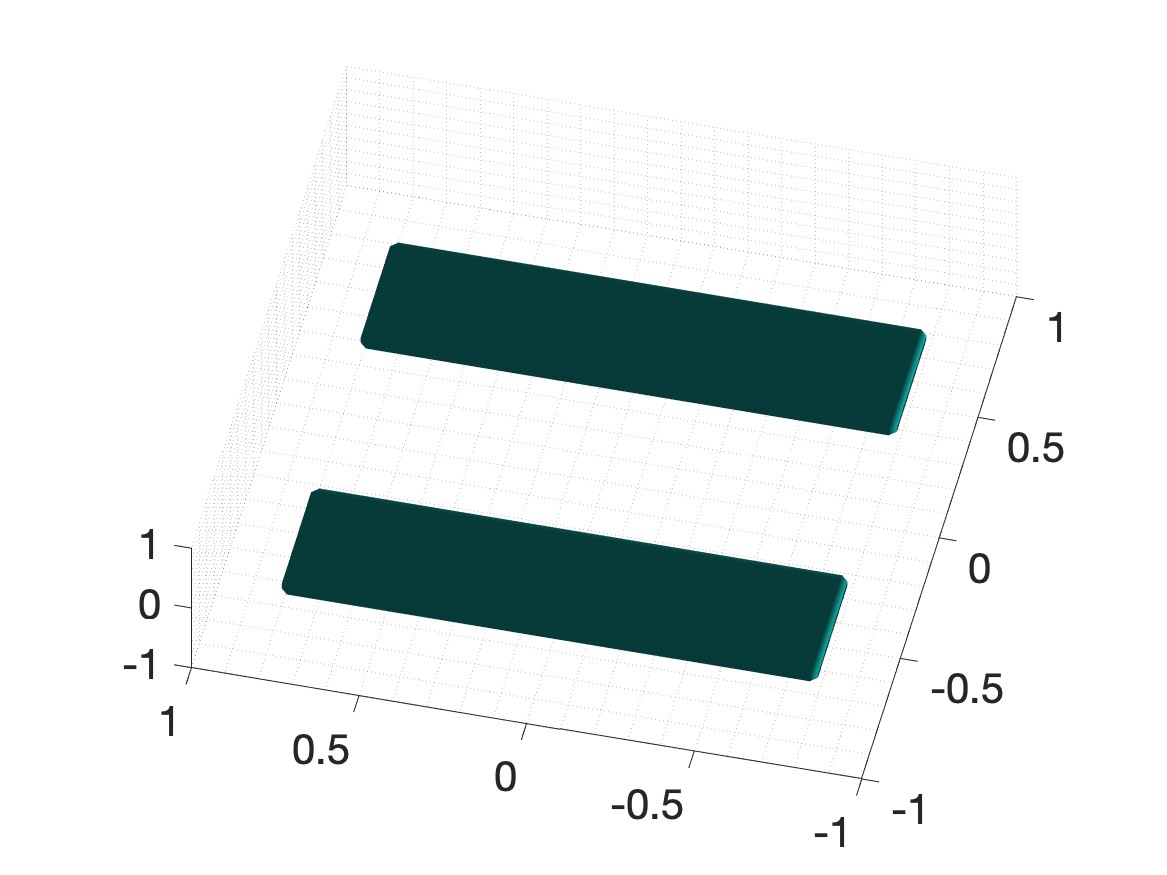}}
	\quad
	\subfloat[\label{test3b}3D view of the reconstructed scatterer]{\includegraphics[width=.45\textwidth]{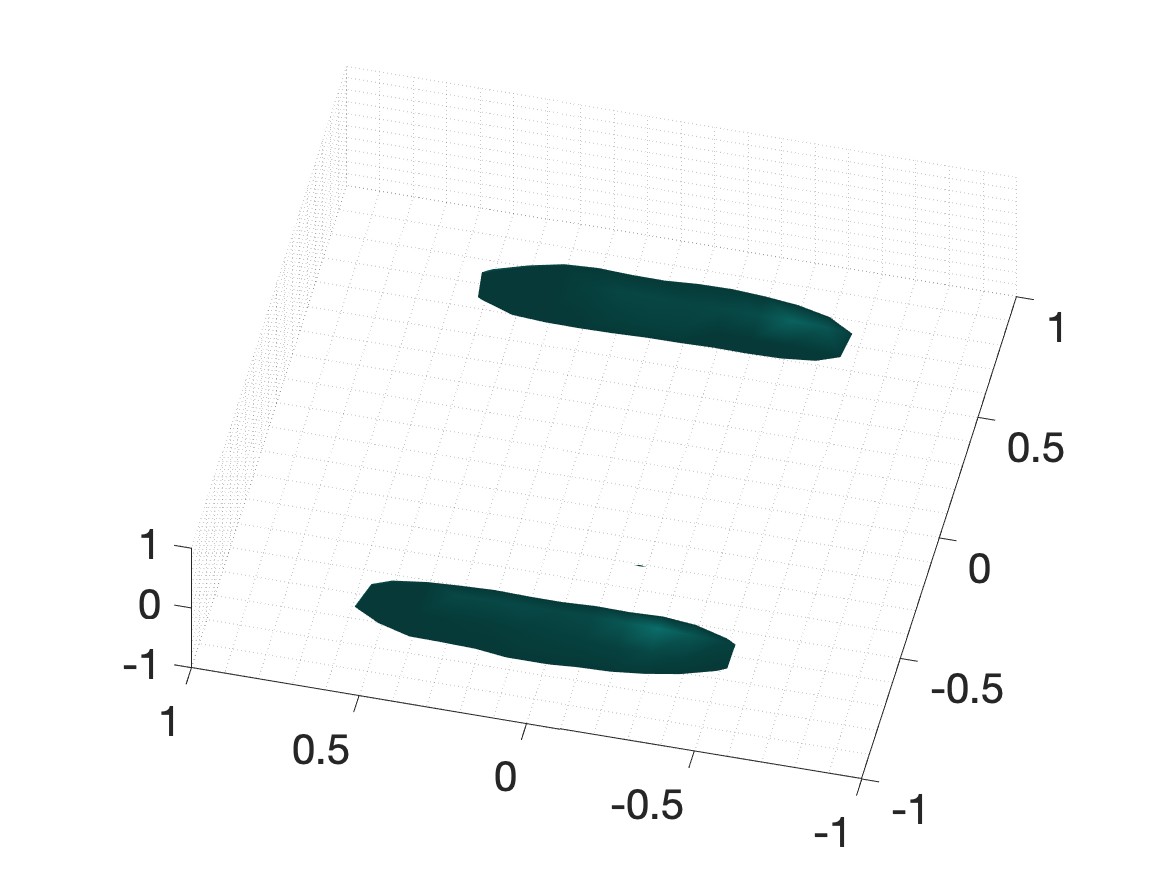}}

	\subfloat[\label{test3c}Cross-sectional view of the function \( c^{\rm true} \)]{\includegraphics[width=.45\textwidth]{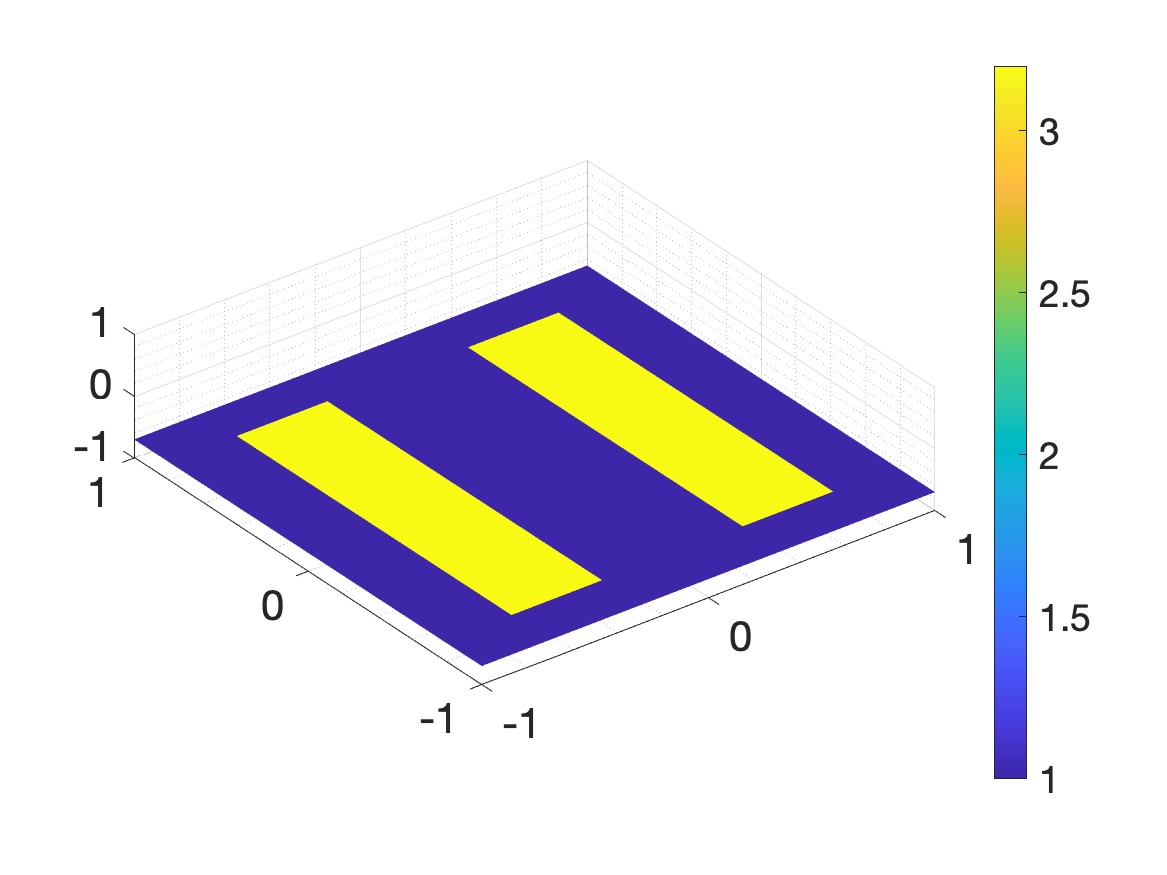}}	
	\quad
	\subfloat[\label{test3d}Cross-sectional view of the function \( c^{\rm comp} \)]{\includegraphics[width=.45\textwidth]{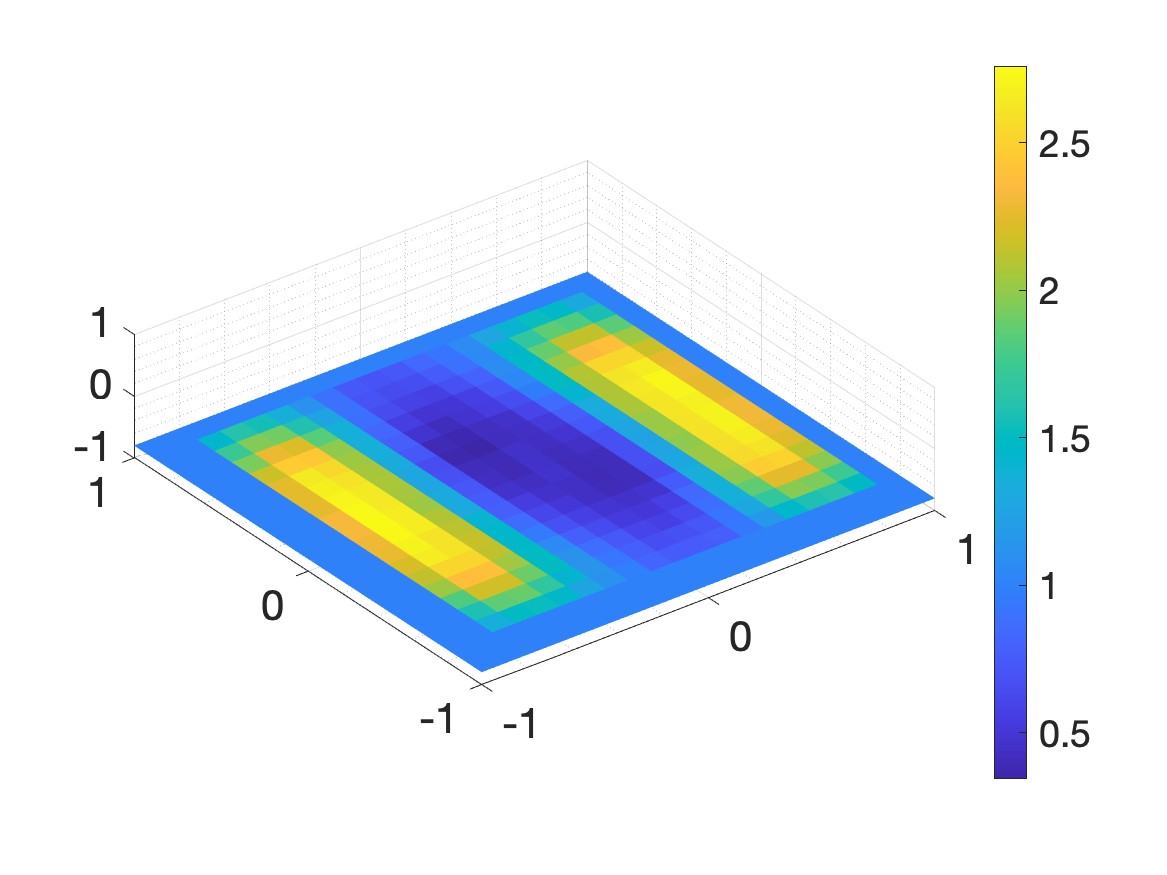}}

	\caption{\label{test3}Test 3. The true and reconstructed dielectric constant \( c(\x) \), \( \x \in \Omega \). Subfigures (a) and (b) show 3D isosurfaces of the true and reconstructed scatterers, respectively, while (c) and (d) display cross-sectional slices of their corresponding dielectric profiles. The reconstruction accurately recovers the elongated shape, spatial separation, and amplitude of both inclusions. These results highlight the method's robustness and effectiveness in handling phaseless data with 10\% noise.}

\end{figure}

Next, we consider the case of two rectangular inclusions. The true dielectric constant function is given by
\[
	c^{\rm true} = 
	\left\{
		\begin{array}{ll}
			3.2 &\max\{0.25|x|, |y + 0.5|\} < 0.2 \mbox{and} |z + 0.65| < 0.05,\\
			3.2 &\max\{0.25|x|, |y - 0.5|\} < 0.2 \mbox{and} |z + 0.65| < 0.05,\\
			1 &\mbox{otherwise}.
		\end{array}
	\right.
\]
The true and reconstructed dielectric profiles are visualized in Figure~\ref{test3}.

Figure~\ref{test3} demonstrates the capability of the proposed method to accurately reconstruct elongated inclusions from noisy phaseless data. The 3D views in subfigures~\ref{test3a} and~\ref{test3b} reveal that the reconstructed scatterers closely resemble the true targets in both shape and orientation, effectively capturing their elongated geometry and spatial arrangement. The cross-sectional slices in subfigures~\ref{test3c} and~\ref{test3d} confirm that the reconstructed dielectric profile approximates the correct locations and amplitudes of the inclusions. 
The maximum value of the computed function $c$ is 2.76, corresponding to the relative noise 13.75\%.
Although some smoothing is visible due to the regularization and the presence of 10\% noise, the method successfully preserves the key features of the underlying structure. These results further validate the robustness of the algorithm, particularly in recovering complex shapes under limited and noisy measurement conditions.

\section{Concluding Remarks}
\label{sec_concl}

In this paper, we have developed a comprehensive numerical framework for solving a 3D phaseless coefficient inverse problem governed by the Helmholtz equation. The method is motivated by both practical imaging applications and a longstanding open question in inverse scattering theory concerning the absence of phase information. Our approach combines several key components: a phase retrieval procedure based on the WKB ansatz, a frequency dimension reduction via truncated Fourier expansion, and the application of the Carleman convexification method to stably reconstruct the spatially varying dielectric constant.

Theoretical results guarantee the strict convexity of the proposed cost functional, and the associated gradient descent method is shown to globally converge to the true solution. Through multiple numerical experiments using simulated noisy data, we have demonstrated that the method accurately recovers both the geometry and contrast of embedded scatterers, even in the presence of high noise and under the challenging constraint of single-sided measurements.

Future research directions include extending this framework to other types of governing equations, such as the full Maxwell system, and further reducing data requirements by exploring compressive sensing or machine learning enhancements. Overall, this work offers a globally convergent and computationally feasible solution to a classically ill-posed and practically relevant inverse problem.

\section*{Acknowledgement}
The work of P.M.N. and L.H.N. was partially supported by the National Science Foundation under grant DMS-2208159.

%
\end{document}